\documentclass[a4paper,twoside,12pt,english]{article}

\usepackage[applemac]{inputenc}
\usepackage[english]{babel}
\makeatletter
\renewcommand\section{\@startsection
{section}{1}{0mm}%
{-2\bigskipamount}%
{\bigskipamount}%
{\normalfont\normalsize\bfseries}%
}

\newcommand\dateymd{\number\year, \ifcase\month\or
January\or February\or March\or April\or May\or June\or
July\or August\or September\or October\or November\or
December\fi, \number\day}
\newcommand\printtime{%
\c@hours=\time \divide\c@hours by60
\c@minutes=\c@hours \multiply\c@minutes by-60
\advance \c@minutes by \time
\ifnum\c@hours<10 0\fi\the\c@hours:%
\ifnum\c@minutes<10 0\fi\the\c@minutes}
\makeatother

\usepackage[OT1]{fontenc}

\usepackage{fancyhdr}
\usepackage{relsize}
\usepackage{bm} 
\usepackage{amssymb}
\usepackage{amsmath}
\usepackage{amsthm}
\usepackage{enumerate}

\usepackage{verbatim} 

\usepackage[nodvipsnames]{color}

\usepackage[pdftex]{graphics} 

\newcommand\upla{}
\newcommand\uplapar{}

\newcommand\bnou{\textup{\textbf{\lower3.7pt\hbox{\char'052}:~}}}
\newcommand\enou{\unskip\textup{\textbf{~:\lower3.7pt\hbox{\char'052}}} }
\newcommand\bvell{\textup{\textbf{\lower3.7pt\hbox{\char'052}:~}$\langle$}}
\newcommand\evell{\unskip\textup{$\rangle$\textbf{~:\lower3.7pt\hbox{\char'052}}} }
\newcommand\bbnou{\textup{\textbf{\lower3.7pt\hbox{\char'052\char'052}:~}}}
\newcommand\eenou{\unskip\textup{\textbf{~:\lower3.7pt\hbox{\char'052\char'052}}} }

\newcommand\ie{i.\,e.~}
\newcommand\ifoi{\,\hbox{if\kern2.5pt and\kern2.5pt only\kern2.5pt if}\,{} }

\newcommand\df{\bfseries}
\newcommand\dfc[1]{\,{\df#1}\,}
\newcommand\dfd[1]{\,{\df#1}\hskip1pt}
\newcommand\secpar[1]{\S\,{#1}}

\newcommand\ensep{\unskip\hskip.65em\ignorespaces}

\newcommand\atilde{\lower3.5pt\hbox{\~{}}}
\newcommand\underl{\lower3.5pt\hbox{-}}

\newcommand\remark{\medskip\noindent\textit{Remark}.\hskip.5em}

\newcommand\better{$\succ$}

\newcommand\pq[2]{\raise.25ex\hbox{\footnotesize${#1}\over{#2}$}%
\hskip-.35ex\null}
\newcommand\onehalf{\leavevmode\raise.5ex\hbox{\scriptsize$1$}\hskip-.25ex/\hskip-.25ex\lower.2ex\hbox{\scriptsize$2$}}

\newcommand\halfsmallskip{\vskip0.5\smallskipamount}

\newcommand\xxxx[1]{%
 \hangindent2.5\parindent
 \hangafter1
 \noindent\hskip.5\parindent
 \hbox to2\parindent{\hss#1\hss}}
\newcommand\condition[2]{\xxxx{#1}\textit{#2}.}
\newcommand\iim[1]{\xxxx{\textup{(#1)}}\ignorespaces}

\newcommand\ddd[1]{\halfsmallskip\vskip-2pt\noindent\hbox to 2\parindent{\hss\footnotesize$\bullet$\ \ }{#1}\ensep}

\newcommand\brwrap[1]{[\textsl{#1}\kern1pt]}

\makeatletter
\newcommand\bibref[1]{\@nameuse{b@#1}}
\renewcommand\@biblabel[1]{\brwrap{#1}}
\renewcommand\@cite[2]{\hbox{\brwrap{#1\if@tempswa\/\upshape\,:\,{\relscale{0.95}#2}\fi}}} 
\makeatother

\newcommand*\dbibref[2]{\bibref{#1}\,\textup{:\,{\relscale{0.95}#2}}}

\newcommand*\refco{\/\kern.1ex\textup{,}\hskip.45ex}
\newcommand*\refsc{\/\kern.15ex\textup{;} }

\newcommand\latop[2]{{#1\atop#2}}

\newcommand\sbset{\subset}

\newcommand\sbseteq{\subseteq}

\newcommand\cd[1]{\!#1\!}

\newcommand\stv{{}^\ast\kern-.25pt v}

\newcommand\hatrho{\hat{%
 \hbox to1.9ex{\hss\rule{0pt}{1ex}\smash{$\rho$}\hskip2pt\hss}}}
\newcommand\hatsigma{\hat{%
 \hbox to1.9ex{\hss\rule{0pt}{1ex}\smash{$\sigma$}\hskip.5pt\hss}}}

\newcommand\tcl[1]{#1{}^\ast}
\newcommand\tclp[1]{\tcl{(#1)}}
\newcommand\crl{\mu}

\newcommand\isc{v^\ast}
\newcommand\icr{\kappa} 
\newcommand\mast{\kappa} 
\newcommand\img{m^\mast}

\newcommand\rxi{\mathrel{\smash{\succ\kern-1.7ex\raise1.15ex\hbox{\mathsurround0pt$\scriptscriptstyle\xi$}\kern.4ex}}}
\newcommand\rxieq{\mathrel{\smash{%
 \vbox{\offinterlineskip\halign{\hfil##\hfil\cr
 \mathsurround0pt$\succ$\cr
 \noalign{\vskip-.5ex}%
 \mathsurround0pt$-$\cr
 \noalign{\vskip-1.15ex}%
 }\vss}\kern-1.05ex\raise1.15ex\hbox{\mathsurround0pt$\scriptscriptstyle\xi$}\kern.4ex}}}

\newcommand\ppmg{m^\sigma}

\newcommand\psc{v^\pi}
\newcommand\pmg{m^\pi}

\newcommand\vk{v^k}

\makeatletter
\gdef\centre#1{\smash{\vbox{\m@th\ialign{\hfil##\hfil\crcr
  \noalign{\kern3\p@}
  \hskip1pt{\tiny$\scriptscriptstyle\bullet$}\crcr
  \noalign{\kern2.5\p@\nointerlineskip}
  $\hfil\displaystyle{#1}\hfil$\crcr}}}}
\makeatother
\newcommand\pcentre[1]{\left(#1\right)^{\raise1pt\hbox{\tiny$\scriptscriptstyle\bullet$}}}

\newcommand\ist{A}
\newcommand\xst{X}
\newcommand\yst{Y}
\newcommand\zst{Z}

\newcommand\tie{\varPi}

\newcommand\istbis{\hbox to1.97ex{\hss\hskip3.5pt$\smash{\widetilde{%
 \hbox to1.9ex{\hss\vphantom{t}\smash{$\ist$}\hskip3.5pt\hss}}}$\hss}}
\newcommand\tiebis{\hbox to1.97ex{\hss\hskip3.5pt$\smash{\widetilde{%
 \hbox to1.9ex{\hss\vphantom{t}\smash{$\tie$}\hskip3.5pt\hss}}}$\hss}}
\newcommand\rhobis{\widetilde\rho}

\newcommand\cst{C}

\newcommand\clustit{\kern.1ex\widetilde c\kern.2ex}
\newcommand\contr[1]{\widetilde #1}

\newcommand\rank[1]{r_{#1}}

\newcommand\arank[1]{\bar r_{#1}}
\newcommand\rlr[1]{R_{#1}}

\newcommand\last{\ell}
\newcommand\anteh{\,{}'\kern-.3ex h}
\newcommand\antef{\,{}'\kern-.45ex f}
\newcommand\antec{\,{}'\kern-.3ex c}

\newcommand\vbis{\widetilde v}

\newcommand\iscbis{\widetilde v{}^{\kern.5pt\ast}}
\newcommand\icrbis{\widetilde\icr}
\newcommand\rlrbis[1]{\hbox to1.97ex{\hss\hskip2.5pt$\smash{\widetilde{%
 \hbox to1.9ex{\hss\vphantom{t}\smash{$R$}\hskip2.5pt\hss}}}$\hss}_{#1}}
\newcommand\pscbis{\widetilde v{}^{\kern.75pt\pi}}

\newcommand\rlrating{\omega}
\newcommand\rlratingbis{\widetilde\omega}

\newcommand\xibis{\smash{\widetilde{\hbox{\vphantom{t}\smash{$\xi$}}}}}
\newcommand\imgbis{\widetilde m^\kappa}
\newcommand\ppmgbis{\widetilde m^\sigma}

\newcommand\gammabis{\widetilde\gamma}

\newcommand\trefbis[1]{\setbox0\hbox{\ref{#1}}%
\smash{\hbox to\wd0{\hss$\widetilde{\hbox to1em{\hss\ref{#1}\hss}}$\hss}}}
\newcommand\pmgbis{\widetilde m^\pi}

\newcommand\gpxy{\gamma\hbox{'}\kern-3pt{}_{xy}}
\newcommand\gbispxy{\widetilde\gamma\hbox{'}\kern-3pt{}_{xy}}

\newcommand\vaa{\mathsf{V}}
\newcommand\vaasub{\vaa\kern-1pt}
\newcommand\vxx{\vaasub_{\scriptscriptstyle X\kern-1pt X}}

\newcommand\vxy{\vaasub_{\scriptscriptstyle X\kern-.25pt Y}}
\newcommand\vyx{\vaasub_{\scriptscriptstyle Y\kern-1pt X}}
\newcommand\vyy{\vaasub_{\scriptscriptstyle Y\kern-.25pt Y}}
\newcommand\vrs{\vaasub_{\scriptscriptstyle R\kern-.25pt S}}
\newcommand\vzz{\vaasub_{\scriptscriptstyle Z\kern-1pt Z}}
\newcommand\vxz{\vaasub_{\scriptscriptstyle X\kern-1pt Z}}

\newcommand\fxx{F_{\scriptscriptstyle X\kern-1pt X}}
\newcommand\fyy{F_{\scriptscriptstyle Y\kern-.25pt Y}}
\newcommand\fxy{F_{\scriptscriptstyle X\kern-.25pt Y}}
\newcommand\fzz{F_{\scriptscriptstyle Z\kern-1pt Z}}
\newcommand\frr{F_{\scriptscriptstyle R\kern-.25pt R}}
\newcommand\flrbis{\hbox to1.75ex{\hss\hskip1.25pt$\widetilde{\hbox to1.5ex{\hss$\varphi$\hskip1.25pt\hss}}$\hss}}
\newcommand\vaabis{\smash{\widetilde{\hbox{\vphantom{t}\smash{$\vaa$}}}}}
\newcommand\vaabissub{\vaabis\kern-1pt}
\newcommand\vaabisxx{\vaabissub_{\scriptscriptstyle X\kern-1pt X}}
\newcommand\vaabiszz{\vaabissub_{\scriptscriptstyle Z\kern-1pt Z}}

\newcommand\xstbiss{\smash{\widetilde\xst}}
\newcommand\ystbiss{\smash{\widetilde\yst}}
\newcommand\zstbiss{\smash{\widetilde\zst}}
\newcommand\vxxbis{\vaabissub_{\scriptscriptstyle \xstbiss\kern-1pt\xstbiss}}
\newcommand\vyxbis{\vaabissub_{\scriptscriptstyle \ystbiss\kern-1pt\xstbiss}}
\newcommand\vzzbis{\vaabissub_{\scriptscriptstyle \zstbiss\kern-1pt\zstbiss}}
\newcommand\flrbissub{\flrbis\kern-1pt}

\newcommand\fbis{\hbox to1.97ex{\hss\hskip2.5pt$\smash{\widetilde{%
 \hbox to1.9ex{\hss\vphantom{t}\smash{$F$}\hskip2.5pt\hss}}}$\hss}}
\newcommand\fbisxx{\fbis_{\kern-2pt\scriptscriptstyle \xstbiss\kern-1pt\xstbiss}}

\newcommand\vbisbis{\smash{\widetilde{\hbox{\vphantom{t}\smash{$\vbis$}}}}\kern.5pt}
\newcommand\iscbisbis{\smash{\widetilde{\hbox{\vphantom{t}\smash{$\vbis$}}\kern.75pt}}^\ast}
\newcommand\icrbisbis{\smash{\widetilde{\hbox{\vphantom{t}\smash{$\icrbis$}}}}\kern.5pt}

\newcommand\chrela{\mathrel{\hbox{\relscale{1.25}$\kern1pt\triangleright\kern1pt$}}}

\newcommand\chrelabis{\mathrel{\smash{\widetilde{\hbox{\vrule width0pt height6.5pt\smash{\hbox{\relscale{1.25}$\kern1pt\triangleright\kern1pt$}}}}}}}

\newcommand\funbis{{\cal F}^{\kern.5pt\prime}}

\newtheorem{proposition}{Proposition}[section]
\newtheorem{lemma}[proposition]{Lemma}
\newtheorem{theorem}[proposition]{Theorem}
\newtheorem{corollary}[proposition]{Corollary}
\newtheorem{open}{Open question} 

\newlength\repskip 
\newenvironment{repeated}[1]%
{\null\hskip-\repskip\hbox to0.9\hsize\bgroup\hbox to0pt{\small$(\ref{#1})$\hss}\hfil\hskip.1\hsize$\displaystyle}%
{$\hfil\egroup}

{\hbox to\hsize\bgroup\hbox to0pt{\small$(\ref{#1})$\hss}\hfil$\displaystyle}%
{$\hfil\egroup\nonumber}

\begin{document}


\thispagestyle{empty}

\null\vskip-5.5\baselineskip\null 
\renewcommand\upla{\vskip-4pt}

\begin{center}
\hrule
\upla
\vskip8mm
\textbf{\uppercase{A continuous rating method for~preferential~voting.\ \ The complete case}}
\par\medskip
\textsc{Rosa Camps,\, Xavier Mora \textup{and} Laia Saumell}
\par
Departament de Matem\`{a}tiques,
Universitat Aut\`onoma de Barcelona,
Catalonia,
Spain
\par\medskip
\texttt{xmora\,@\,mat.uab.cat}
\par\medskip
July 13, 2009;\ensep revised February 24, 2011 
\vskip7.5mm
\upla
\hrule
\end{center}

\upla
\vskip-4mm\null
\begin{abstract}
A method is given for quantitatively rating the social acceptance of different options which are the matter of a complete preferential
\linebreak[3]
vote. Completeness means that every voter expresses a~comparison 
(a~preference or a~tie) about each pair of options.
The~proposed method is proved to have certain desirable properties, which include:\ensep
the conti\-nuity of the rates with respect to the data,\ensep
a decomposition property that characterizes certain situations opposite to a tie,\ensep
the Condorcet-Smith principle,\ensep 
and
clone consistency.\ensep
One can view this rating method as a complement for the ranking method introduced in~1997 by Markus Schulze.
It~is also related to certain methods of one-dimensional scaling or cluster analysis.

\upla

\upla
\bigskip\noindent
\textbf{AMS subject classifications:} 
\textit{%
05C20, 
91B12, 
91B14, 
91C15, 
91C20. 
}
\end{abstract}

\vskip5mm
\upla
\hrule

\vskip-12mm\null 
\section*{}
The outcome of a vote is commonly expected
to specify not only a winner and an ordering of the candidates,
but also a quantitative estimate of the social acceptance of each of them.
Such a quantification is expected even 
when the individual votes give only qualitative information.
%

The simplest voting methods are clearly based upon such a quantification.
This is indeed the case of the plurality count as well as that of the Borda count.
\ensep
However, it is well known that
these methods do not comply with basic majority principles
nor with other desirable conditions.
\ensep
In~order~to satisfy certain combinations of such principles and conditions,
one~must resort to other more elaborate methods, such~as 
the celebrated rule of Condorcet, Kem\'eny and Young
\hbox{\brwrap{
\bibref{mo}\refsc
\dbibref{t6}{p.\,182--190}%
}},
the method of ranked pairs
\hbox{\brwrap{
\bibref{ti}\refco
\bibref{zt}\refsc
\dbibref{t6}{p.\,219--223}%
}},
or the method introduced in~1997 by~Markus Schulze
\hbox{\brwrap{
\bibref{sc}\refco
\bibref{scbis}\refsc
\dbibref{t6}{p.\,228--232}%
}},
which we will refer to as the method of paths.
\ensep
Now,~as they stand, 
these methods rank the candidates in a purely ordinal way,
without properly quantifying the social acceptance of each of them.
\ensep
%
So,~it is natural to ask for a method that
combines the above-mentioned
principles and conditions 
with a quantitative rating of the candidates.
%
%

\smallskip
A quantitative rating should allow to sense the closeness between two candidates,
such as the winner and the runner-up.
For that purpose, it is essential that the rates vary in a \emph{continuous} way,
especially through situations where ties or multiple orders occur.

On~the other hand,
it should also allow to recognise
certain 
situations that are opposite to a~tie.
For instance, a candidate should get the best possible rate
\ifoi it has been placed first by all voters.
This is a particular case of a more general condition that we will call \emph{decomposition}. This condition, that will be made precise later on, places sharp constraints on the rates that should be obtained when the candidates are partitioned in two classes~$\xst$ and~$\yst$ such that each member of $\xst$ is unanimously preferred to every member of~$\yst$.

\newcommand\llmp{\textup{M}}

In this article we will produce a rating method that combines such a quantitative character with other desirable properties of a qualitative nature.
Among them we will be especially interested in the following extension of the Condorcet principle introduced in 1973 by John H.~Smith \cite{smith}:
Assume that the set of candidates is partitioned in two classes~$\xst$ and~$\yst$ such that for each member of $\xst$ and every member of~$\yst$ there are more than half of the individual votes where the former is preferred to the latter; in that case,
the social ranking should also prefer each member of $\xst$ to any member of $\yst$.
%
%
This principle is quite pertinent when one is interested not only in choosing a winner but also in ranking all the alternatives (or in rating them).

To our knowledge, the existing literature does not offer any other rating method that combines this principle with the above-mentioned quantitative properties of continuity and decomposition.
\ensep
We will refer to our method as the \dfd{CLC~rating method}, where the capital letters stand for ``Continuous Llull Condorcet''.
\ensep
The reader interested to try it can use the \,\textsl{CLC~calculator}\, which has been made available at~\cite{clc}.

\smallskip
Of course, any rating automatically implies a ranking.
In this connection, it should be noticed that
the CLC~rating method is built upon Schulze's method of paths
as underlying ranking method.
As we will remark in the concluding section, we doubt that any of the other ranking methods mentioned above
could be extended
to a rating method with the properties of continuity and decomposition.
\ensep
Having said that, it should be clear 
that the present work is not aimed at saying anything new about ranking methods as such.
Whatever we might say about them will always be in reference to the rating issue.


\smallskip
By reasons of space, this article is restricted to the complete case and to a particular class of rates that we call rank-like rates. We~are in the complete case when every individual expresses a comparison (a preference or a tie) about each pair of options. The incomplete case requires some additional developments that are dealt with in a separate article~\cite{cri}.
In~another separate article we deal with another class of rates that have a fraction-like character~\cite{crz}.

\smallskip
The present article is organized as follows: 
Section~1 gives a more precise statement of the problem and finishes with a general remark. 
Section~2 presents an heuristic outline of our proposal, ending with a summary of the procedure and an illustrative example. 
Section~3 introduces some mathematical language. 
Sections~\ref{sec-indirect-scores}--\ref{sec-condorcet} give detailed mathematical proofs of the claimed properties.\ensep
Sections~\ref{sec-clons}--\ref{sec-monotonicity} are devoted to other interesting properties of the concomitant social ranking, namely clone consistency and two weak forms of monotonicity.\ensep
Finally, section~\ref{sec-conclusion} makes some concluding remarks and poses a few open questions.

\section{Statement of the problem and a general remark}

\paragraph{1.1}
Let us consider a set of $N$~options which are the matter of a vote.
Let us assume that each voter expresses his preferences in a qualitative way,
for instance by listing those options in order of preference.
Our aim is to combine such individual preferences so as to rate the social acceptance of each option on a continuous scale.
More specifically, we would like to do it in accordance with
the following conditions:

\smallskip
\newcommand\llsi{\textup{A}}
\condition{\llsi}{Scale invariance (homogeneity\hskip.1pt)}
The rates depend only on the relative frequency of each possible content of an individual vote.
In~other words, if every individual vote is replaced by a fixed number of copies of it, the rates remain exactly the same.

\smallskip
\newcommand\llpe{\textup{B}}
\condition{\llpe}{Permutation equivariance (neutrality\hskip.1pt)}
Applying a certain permutation of the options to all of the individual votes has no other effect than getting the same permutation in the social rating.

\smallskip
\newcommand\llco{\textup{C}}
\condition{\llco}{Continuity}
The rates depend continuously on the relative frequency of each possible content of an individual vote.

\newcommand\condr[1]{#1\rlap{\mathsurround0pt$_r$}}
\newcommand\condf[1]{#1\rlap{\mathsurround0pt$_f$}}

\pagebreak 
\null\vskip-14mm\null

\medskip
The next two conditions choose a specific form of rating. From now on we will refer to it as rank-like rating. In the complete case considered in the present paper, these two conditions take the following form:

\smallskip
\newcommand\llrr{\textup{D}}
\condition{\llrr}{Rank-like form (complete case)}
Each rank-like rate is a number, integer or fractional, between $1$ and~$N$. The~best possible value is~$1$ and the worst possible one is~$N$. The average rank-like rate is~$(N+1)/2$.

\smallskip
\newcommand\llrd{\textup{E}}
\condition{\llrd}{Rank-like decomposition (complete case)}
Consider a splitting of the options in two classes~$\xst$ and~$\yst$. Consider the case where each member of $\xst$ is unanimously preferred to every member of~$\yst$.
This fact is equivalent to each of the following ones,
where $|\xst|$ denotes the number of elements of $\xst$:
\ensep
(a)~The rank-like rates of~$\xst$ coincide with those that one obtains when the individual votes are restricted to~$\xst$.\linebreak 
\ensep
(b)~After diminishing them by the number~$|\xst|$, the rank-like rates of~$\yst$ coincide with those that one obtains when the individual\linebreak 
votes are restricted to~$\yst$.
\ensep
(c)~The average rank-like rate of~$\xst$ is $(|\xst|+1)/2$.

\smallskip
\noindent
In particular, an option will get a rank-like rate exactly equal to~$1$ \,if and only if\, it is unanimously preferred to any other.
Similarly, an option will get a rank-like rate exactly equal to~$N$\,if and only if\, it is unanimously considered worse than any other.

\medskip
Finally, we require a condition that concerns only
the concomitant social ranking, that is,
the ordinal information contained in the social rating:

\smallskip
\condition{\llmp}{Condorcet-Smith principle} 
Consider a splitting of the options in two classes~$\xst$ and~$\yst$. 
Assume that for each member of $\xst$ and every member of $\yst$ there are more than half of the individual votes where the former is preferred to the latter.
In that case, the social ranking also prefers each member of $\xst$ to every member of~$\yst$.

\renewcommand\upla{\vskip-6pt}
\upla
\paragraph{1.2}
Let us emphasize that the individual votes that we are dealing with do \textit{not} have a quantitative character (at least for the moment): each voter is allowed to express a preference for $x$ rather than $y$, or vice versa, or maybe a tie between them, but he is not allowed to quantify such a preference.

This contrasts with ``range voting'' methods, where each individual vote is already a quantitative rating
\hbox{\brwrap{
\dbibref{t6}{p.\,174--176}\refsc
\bibref{bala}%
}}.
Such methods are free from many of the difficulties that lurk behind the present setting.
However, they make sense only as long as all voters mean the same by  each possible value of the rating variable.
This hypothesis may be reasonable in some cases, but in many others it is hardly valid.
In fact, voting is often used in connection with moral, psychological or aesthetic qualities, whose appreciation may be as little quantifiable, but also as much ``comparable'', as, for instance, the feelings of pleasure or pain.

In our case, the quantitative character of the output will derive from the fact of having a certain number of qualitative preferences.
The~larger this number, the more meaningful will be the quantitative character of the social rating.
This is especially applicable to the continuity property~\llco,
according to which
a small variation in the proportion of votes with a given content produces only small variations in the rates.
In fact, a few votes will be a small proportion only in the measure that 
the total number of votes is large enough.



\section{Heuristic outline}\label{sec-heuristic}

This section presents our proposal as the result of a quest
for the desired properties. Hopefully, this will communicate the main ideas that lie behind the formulas.

\paragraph{2.1}
The aim of complying with condition~\llmp\ 
calls for the point of view of \dfd{paired comparisons}.
So our starting point will be the numbers $V_{xy}$ that count how many voters prefer $x$ to $y$. In this connection we will take the view that each vote that ties $x$ with $y$ is equivalent to half a vote preferring $x$ to $y$\, plus another half a vote preferring~$y$~to~$x$.
In order to achieve scale invariance, we will immediately switch to the corresponding fractions $v_{xy} = V_{xy}/V$, where $V$ denotes the total number of votes.
\ensep
The numbers~$v_{xy}$ will be called the binary \dfc{scores} of the vote, and their collection will be called the \dfc{Llull matrix} of the vote.
\ensep
Since we are considering the case of complete votes, these numbers are assumed 
to satisfy
\begin{equation}
v_{xy} + v_{yx} \,=\, 1.
\label{eq:complete}
\end{equation}

\medskip
Besides the scores $v_{xy}$, in the sequel we will often deal with the
\dfc{margins} $m_{xy}$, which are defined by
\begin{equation}
m_{xy} \,=\, v_{xy} - v_{yx}.
\end{equation}
Obviously, their dependence on the pair $xy$ is antisymmetric, that is
\begin{equation}
m_{yx} \,=\, - m_{xy}.
\end{equation}
It is clear also that the equality~(\ref{eq:complete}) allows to recover the scores from the margins by means of the following formula:
\begin{equation}
v_{xy} \,=\, (1 + m_{xy})/2.
\end{equation}

\paragraph{2.2}
A natural candidate for defining the social preference is the following: $x$~is socially preferred to~$y$
whenever $v_{xy} > v_{yx}$.
Of course, it can happen that $v_{xy} = v_{yx}$,
in~which case one would consider that $x$~is socially tied with~$y$. 
The~binary relation that includes all pairs $xy$ for which $v_{xy} >
v_{yx}$ will be denoted by $\crl(v)$ and will be called the
\dfd{comparison relation};
together with it, we will consider also the relation 
$\hat\crl(v)$ defined by the non-strict inequality $v_{xy} \ge v_{yx}$.

\medskip
As it is well known, the main problem with paired comparisons
is that the comparison relations $\crl(v)$ and~$\hat\crl(v)$
may lack transitivity even if the individual
preferences are all of them transitive~\cite{mu,t6}.

\paragraph{2.3}
The next developments rely upon an operation $(v_{xy})\rightarrow (\isc_{xy})$
that transforms the original system of binary scores into a new one.
This operation is defined in the following way:
for every pair~$xy$, one considers all possible paths $x_0 x_1 \dots x_n$
going from~$x_0 = x$ to~$x_n = y$; every such path is associated with the
score of its weakest link, \ie the smallest value of $v_{x_ix_{i+1}}$;
finally, $\isc_{xy}$ is defined as the maximum value of this associated score
over all paths from $x$ to $y$.
In other words,
\begin{equation}
\isc_{xy} \hskip.75em = \hskip.75em
\max_{\vtop{\scriptsize\halign{\hfil#\hfil\cr\noalign{\vskip.5pt}$x_0=x$\cr$x_n=y$\cr}}}
\hskip.75em
\min_{\vtop{\scriptsize\halign{\hfil#\hfil\cr\noalign{\vskip-1.25pt}$i\ge0$\cr$i<n$\cr}}}
\hskip.75em v_{x_ix_{i+1}},
\label{eq:paths}
\end{equation}
where the \,$\max$\, operator considers all possible paths from $x$ to $y$,
and the \,$\min$\, operator considers all the links of a particular path.
The scores $\isc_{xy}$ will be called the \dfc{indirect scores}
associated with the (direct) scores $v_{xy}$.

If $(v_{xy})$ is the table of 0's and 1's associated with a binary relation $\rho$ (by putting $v_{xy}=1$ \ifoi $xy\in\rho$), then $(\isc_{xy})$
is exactly the table associated with $\rho^\ast$,
the transitive closure of $\rho$.
So, the operation $(v_{xy})\mapsto (\isc_{xy})$ can be viewed as a
quantitative analogue of the notion of transitive
closure (see~\cite[Ch.\,25]{co}).

The main point, remarked in 1998 by Markus Schulze~%
\hbox{\brwrap{\bibref{sc}\,b}}, 
is that the comparison relation
associated with a table of indirect scores is always transitive (Theorem~\ref{st:transSchulze}).
So, $\crl(\isc)$ is always transitive,
no matter what the case is for~$\crl(v)$.
This is true in spite of the fact that
$\crl(\isc)$ can easily differ from~$\crl^\ast(v)$ (the transitive closure of $\crl(v)$).

\paragraph{2.4}
In the following we put
\begin{equation}
\icr=\crl(\isc),\qquad \hat\icr=\hat\crl(\isc),\qquad
\img_{xy} = \isc_{xy} - \isc_{yx}.
\end{equation}
So, $xy\in\icr$ \ifoi $\isc_{xy} > \isc_{yx}$, \ie $\img_{xy}>0$,
and $xy\in\hat\icr$ \ifoi $\isc_{xy} \ge \isc_{yx}$, \ie $\img_{xy}\ge0$.
From now on we will refer to~$\icr$ as the \dfd{indirect comparison relation},
and to $\img_{xy}$ as the {\df indirect margin} associated with the pair $xy$.

As it has been stated above,
the relation~$\icr$ is transitive.\ensep
Besides that, it~is clearly asymmetric (one cannot have
both $\isc_{xy} > \isc_{yx}$ and vice versa).
On~the other hand, it may be incomplete
(one can have $\isc_{xy} = \isc_{yx}$).
When it differs from~$\icr$, the complete relation $\hat\icr$
is not asymmetric and ---somewhat surprisingly---
it may be not transitive either.
\ensep
However, one can always find a total order~$\xi$ that
satisfies $\icr \sbseteq \xi \sbseteq \hat\icr$ (Theorem~\ref{st:existenceXiThm}).
From now on, any total order~$\xi$ that satisfies this condition
will be called an \dfc{admissible order}.
Let us remark that such a definition is redundant: in fact, one easily sees that each of the required inclusions implies the other one.

The rating that we are looking for will be based on such an order~$\xi$.
More specifically, it will be compatible with $\xi$ in the sense that the rank-like rates~$\rlr{x}$
will satisfy the inequality $\rlr{x}\le\rlr{y}$ whenever $xy\in\xi$.
If $\icr$ is already a total order,
so that $\xi=\icr$,
the preceding inequality will be satisfied in the strict form $\rlr{x}<\rlr{y}$,
and this will happen \ifoi $xy\in\icr$.


The following steps assume that one has fixed an admissible order~$\xi$.
From now on the situation~$xy\in\xi$ will be expressed also by~$x\rxi y$.
According to the definitions, the inclusions $\icr \cd\sbseteq \xi \cd\sbseteq \hat\icr$
are~equivalent to saying that $\isc_{xy}>\isc_{yx}$ implies $x\rxi y$
and that the latter implies $\isc_{xy} \ge \isc_{yx}$.
In other words, if the different options are ordered according to $\rxi$,
the~matrix $\isc_{xy}$ has then the property that each element above the diagonal
is larger than or equal to its symmetric over the diagonal.

\paragraph{2.5}
Rating the different options means positioning them on a line.
Besides\linebreak[3] complying with the qualitative restriction of
being compatible with~$\xi$ in the sense above,
we want that the distances between items
reflect the quantitative information provided by the binary scores.
However, a rating is expressed by $N$~numbers,
whereas the binary scores are $N(N-1)$ numbers.
So we are bound to do some sort of projection.
Problems of this kind have a certain tradition in
combinatorial data analysis and cluster analysis~\cite{js,ham}.
In~fact, some of the operations that will be used below
can be viewed from that point of view.

\medskip
Let us assume for a while that
the votes are total orders,
\ie each vote lists all the options by order of preference, without any ties.
This is the standard case for the application of Borda's method,
which is linearly equivalent to rating each option by the mean value of its ranks,
\ie the~ordinal numbers that give its position in these different orders.
As it was noticed by Borda himself
(in his setting linearly related to ours),
these mean ranks, which we will denote by~$\arank{x}$,
can be obtained from the Llull matrix by means of the following formula:
\begin{equation}
\arank{x} \,=\, N - \sum_{y\neq x}\,v_{xy},
\label{eq:avranksfromscores}
\end{equation}
or equivalently,
\begin{equation}
\arank{x} \,=\, (N+1 - \sum_{y\neq x}\,m_{xy}\,) \,/\, 2.
\label{eq:avranksfrommargins}
\end{equation}

Let us look at the meaning of the margins~$m_{xy}$ in connection with the idea of projecting
the Llull matrix into a rating: If there are no other items than $x$~and~$y$,
we can certainly view the sign and magnitude of~$m_{xy}$ as giving respectively the qualitative and quantitative aspects of the relative positions of $x$ and $y$ on the rating line,
that~is, the order and the distance between them.
\ensep
When there are more than two items, however,
we have several pieces of information of this kind, one for every pair,
and these different pieces may be incompatible with each other,
quantitatively or even qualitatively,
which motivates indeed the problem that we are dealing with.
In particular, the~mean ranks~$\arank{x}$ often violate
the desired compatibility with the relation~$\xi$.

\medskip
In order to construct a rating compatible with~$\xi$, we will
use a formula analogous to~(\ref{eq:avranksfromscores}) where
the scores~$v_{xy}$ are replaced by certain \dfc{projected scores} $\psc_{xy}$ to be defined in the following paragraphs.
Together with them, we will make use of the corresponding \dfc{projected margins}~$\pmg_{xy} = \psc_{xy} - \psc_{yx}$.
Like the original scores (but not necessarily the indirect ones) the projected scores will be required to satisfy the equality $\psc_{xy}+\psc_{yx} = 1$, from which it follows that $\psc_{xy} = (1+\pmg_{xy})/2$.
\ensep
So,~the 
rates that we are looking for will be obtained in the following way:
\begin{equation}
\rlr{x} \,=\, N - \sum_{y\neq x}\,\psc_{xy},
\label{eq:rrates}
\end{equation}
or equivalently,
\begin{equation}
\rlr{x} \,=\, (N+1 - \sum_{y\neq x}\,\pmg_{xy}\,) \,/\, 2.
\label{eq:rratesfrommargins}
\end{equation}
Such formulas will be used not only in the case where the votes are total orders, but also in more general situations.

\paragraph{2.6}
Our goals will be achieved by defining the projected margins in the following way, where we assume $x\rxi y$ and $x'$~denotes the item that immediately follows~$x$ in the total order~$\xi$:
\begin{gather}
\img_{xy} \,=\, \isc_{xy} - \isc_{yx},
\label{eq:cprojection1}
\\[2.5pt]
\ppmg_{xy} \,=\, \min\,\{\, \img_{pq} \;\vert\; p\rxieq x,\; y\rxieq q\,\},
\label{eq:cprojection2}
\\[2.5pt]
\pmg_{xy} \,=\, \max\,\{\, \ppmg_{pp'} \;\vert\; x\rxieq p\rxi y\,\},
\label{eq:cprojection3}
\\[2.5pt]
\pmg_{yx} \,=\, -\pmg_{xy}.
\label{eq:cprojection4}
\end{gather}

\medskip
As one can easily check, this construction ensures that
\begin{equation}
\pmg_{xz} \,=\, \max\,(\pmg_{xy},\pmg_{yz}),\qquad \hbox{whenever $x\rxi y\rxi z$}.
\label{eq:equaltomax}
\end{equation}
From this equality it follows that the absolute values $d_{xy} = |\pmg_{xy}|$ satisfy the following inequality,
which makes no reference to the relation $\xi$:
\begin{equation}
d_{xz} \,\le\, \max\,(d_{xy},d_{yz}),\qquad \hbox{for any $x,y,z$.}
\label{eq:ultrametric}
\end{equation}
This condition, called the ultrametric inequality,
is well known in cluster analysis,
where it appears as a necessary and sufficient condition
for the dissimilarities $d_{xy}$ to define a 
hierarchical classification of the set under consideration
\hbox{\brwrap{
\dbibref{js}{\secpar{7.2}}\refsc
\dbibref{ham}{\secpar{3.2.1}}%
}}.

\remark
The operation $(\img_{xy})\rightarrow(\pmg_{xy})$ defined by (\ref{eq:cprojection2}--\ref{eq:cprojection3}) is akin to the single-link method of cluster analysis, which can be viewed as a continuous method for projecting a matrix of dissimilarities onto the set of ultrametric distances; such a continuous projection is achieved by taking the maximal ultrametric distance which is bounded by the given matrix of dissimilarities \cite[\secpar{7.3,\,7.4,\,8.3,\,9.3}]{js}.
The operation $(\img_{xy})\rightarrow(\pmg_{xy})$ does the same kind of job under the constraint that the clusters 
be intervals of the total order~$\xi$.

\paragraph{2.7}
\textbf{Summary of the procedure}
\begin{enumerate}
\setlength\itemsep{0pt}
\setcounter{enumi}{-1}
\item Form the Llull matrix~$(v_{xy})$~(\secpar{2.1}).
\item Compute the indirect scores~$\isc_{xy}$ defined by~(\ref{eq:paths}). An efficient way to do it is the Floyd-Warshall algorithm \cite[\secpar{25.2}]{co}. Work out the indirect margins $\img_{xy} \cd= \isc_{xy}\cd-\isc_{yx}$.
\item Consider the indirect comparison relation $\icr = \{xy\mid\img_{xy} > 0\}$.\linebreak
Fix an admissible order~$\xi$, \ie a total order that extends $\icr$. For instance, it suffices to arrange the options by non-decreasing values of
the ``tie-splitting'' Copeland scores
$\rank{x} = N - |\{\,y\mid y\cd\neq x,\ \img_{xy}\cd>0\}|\linebreak - \frac12\,|\{\,y\mid y\cd\neq x,\ \img_{xy}\cd=0\}|$ (Proposition~\ref{st:Copeland}).
\item Starting from the indirect margins~$\img_{xy}$, work~out the superdiagonal intermediate projected margins $\ppmg_{xx'}$ as defined in (\ref{eq:cprojection2}).
\item Compute the projected margins $\pmg_{xy}$ according to~(\ref{eq:cprojection3}--\ref{eq:cprojection4}). The projected scores are then determined by the formula $\psc_{xy}=(1+\pmg_{xy})/2$.
\item Compute the rank-like rates~$\rlr{x}$ according to (\ref{eq:rrates}) (here equivalent to~(\ref{eq:rratesfrommargins})).
\end{enumerate}

\smallskip
The computing time is of order $N^3$, where $N$ is the number of options.\ensep 
The \textsl{CLC calculator} made available at \cite{clc} allows to follow the details of the procedure by choosing the option ``Detailed mode''.

\newbox\Strutbox
\newdimen\Strutheight

\def\initsize{%
 \Strutheight=.7\baselineskip \advance\Strutheight by2pt
 \setbox\Strutbox=\hbox{\vrule height\Strutheight depth .3\baselineskip width0pt}%
 \def\Strut{\relax\unhcopy\Strutbox}%
}
\initsize

\def\initsizesmall{%
 \Strutheight=.7\baselineskip \advance\Strutheight by2pt
 \setbox\Strutbox=\hbox{\vrule height\Strutheight depth .3\baselineskip width0pt}%
 \def\Strut{\relax\unhcopy\Strutbox}%
}

\edef\bv{|} \catcode`\|=\active \let|\bv
\edef\bi{/} \catcode`\/=\active \let/\bi
\def\htskip{\hskip.4em plus.25em minus.25em}
\newcount\nspan

\def\ls#1{\hbox{{\xipt\sl@\/}#1}}
\def\Ls#1{\hbox{{\xipt\sl@\/}#1}}
\def≈{\,=\,}

\newif\ifhborders

\def\taula#1{
\vbox\bgroup
\def|{&&}
\ifhborders
\def/{\unskip&\cr\noalign{\hrule height.2pt}\Strut&&\ignorespaces}
\else
\def/{\unskip&\cr\strut&&\ignorespaces}
\fi
\def\hrf{\unskip&\cr\noalign{\hrule height.2pt}\Strut&&\ignorespaces}
\def\hrg{\unskip&\cr\noalign{\hrule height.8pt}\Strut&&\ignorespaces}
\halign\bgroup##\hfil\cr 
\hbox{\sc\ #1\hss}\cr 
\noalign{\vskip5pt}
\hbox\bgroup\vrule
}

\def\fitaula{\vrule\egroup\cr\egroup\egroup}

\def\bloc #1#2#3#4#5#6#7\ska{%
\nspan=#1\multiply\nspan by2\advance\nspan by-1%
\vrule width.2pt%
\hbox{\vbox{\offinterlineskip\halign{%
##&\vrule## width.2pt%
&&\htskip#2{#4##\/}#3\htskip &\vrule## width.2pt\cr
\noalign{\hrule height.8pt}
\Strut&&\multispan\nspan\it#5\/\hfil&\cr
\strut&&#6&\cr
\noalign{\hrule height.8pt}
\Strut&&\ignorespaces#7&\cr
\noalign{\hrule height.8pt}
}}}%
\vrule width.2pt}

\def\multicol(#1)#2:#3\ska{\bloc{#1}\hfil\hfil\rm\null{#2}{#3}\ska}
\def\hmulticol(#1)#2:#3:#4\ska{\bloc{#1}\hfil\hfil\rm{#2}{#3}{#4}\ska}
\def\gap#1{\vrule\egroup{#1}\hbox\bgroup\vrule}

\def\numeros#1\ska{\bloc1\hfil\hfil\sf\null{\it Id\/}#1\ska} 
\def\numx#1\ska{\bloc1\hfil\hfil\sf\null{$x$}#1\ska} 
\def\jutges(#1)#2:#3\ska{\bloc{#1}\hfil\hfil\rm{Judges}{#2}{#3}\ska}
\def\calculs(#1)#2:#3\ska{\bloc{#1}\hfil\hfil\rm{Computations\ }{#2}{#3}\ska}
\def\balls(#1)#2:#3\ska{\bloc{#1}\hfil\hfil\rm{Dances}{#2}{#3}\ska}
\def\rxxi(#1)#2:#3\ska{\bloc{#1}\null\hfil\rm{Computations}{#2}{#3}\ska}
\def\balljutges#1(#2)#3:#4\ska{\bloc{#2}\hfil\hfil\rm{#1}{#3}{#4}\ska}
\def\fballs(#1)#2:#3\ska{\bloc{#1}\null\hfil\rm{Dances}{#2}{#3}\ska}
\def\rate#1\ska{\bloc1\hfil\hfil\bf\null{$X$}{#1}\ska}

\def\fl{\it} 
\def\al{\sl} 

\def\llocs#1\ska{\bloc1\hfil\hfil\fl\null{\it R\/}#1\ska}
\def\suma#1\ska{\bloc1\hfil\hfil\bf\null{\it S\/}#1\ska}
\def\promig#1\ska{\bloc1\hfil\hfil\bf\null{\it A\/}{#1}\ska}
\def\reordenats(#1):#2\ska{\bloc{#1}\hfil\hfil\rm\null%
{\multispan\nspan\hfil{\it Rearranged\/}\hskip3.5pt\hfil}{#2}\ska}
\def\mediana#1\ska{\bloc1\hfil\hfil\bf\null{\it M\/}#1\ska}
\def\coco#1\ska{\bloc1\hfil\hfil\bf\null{\it C\/}#1\ska}

\def\dmedas(#1)#2:#3\ska{\bloc{#1}\hfil\hfil\rm{Dances}{#2}{#3}\ska}
\def\jmedas(#1)#2:#3\ska{\bloc{#1}\hfil\hfil\rm{Judges}{#2}{#3}\ska}
\def\bmedas(#1)#2:#3\ska{\bloc{#1}\hfil\hfil\rm{Balance}{#2}{#3}\ska}

\def\mballs(#1)#2:#3\ska{\bloc{#1}\hfil\hfil\rm{Dances \rm($M$)}{#2}{#3}\ska}
\def\md#1\ska{\bloc1\hfil\hfil\bf\null{$M^D$}#1\ska}
\def\mj#1\ska{\bloc1\hfil\hfil\bf\null{$M^J$}#1\ska}
\def\mb#1\ska{\bloc1\hfil\hfil\bf\null{$M^B$}#1\ska}
\def\liballs(#1)#2:#3\ska{\bloc{#1}\hfil\hfil\rm{Dances \rm($L_1$)}{#2}{#3}\ska}
\def\lid#1\ska{\bloc1\hfil\hfil\rm\null{$L_1^D$}#1\ska}

\def\xd#1\ska{\bloc1\null\hfil\bf\null{$X^D$}{#1}\ska}
\def\xj#1\ska{\bloc1\null\hfil\bf\null{$X^J$}#1\ska}
\def\xdj#1\ska{\bloc{2}\hfil\hfil\rm\null{$X^D$|$X^J$}{#1}\ska}
\def\xb#1\ska{\bloc1\null\hfil\bf\null{$X^B$}#1\ska}

\def\rates#1\ska{\bloc1\hfil\hfil\rm\null{\it X\/}#1\ska}
\def\ratesbf#1\ska{\bloc1\hfil\hfil\bf\null{\it X\/}#1\ska}

\def\peu{\unskip&\cr\noalign{\hrule height.8pt}\Strut&&\ignorespaces}
\def\pballs#1:#2\ska{\bloc1\hfil\hfil\rm{\sf #1 }\null
{#2\peu $\scriptstyle D$/$\scriptstyle J$\peu$\scriptstyle B$}\ska}

\def\hbi#1{\hbox to.55em{\hss #1\hss}}
\def\hbij#1{\hbox to.85em{\hss #1\hss}}
\def\hbii#1{\hbox to1em{\hss #1\hss}}
\def\hbiii#1{\hbox to1.3em{\hss #1\hss}}

\def\iiijutges{\hbi{\sf A}|\hbi{\sf B}|\hbi{\sf C}}
\def\ivjutges{\hbi{\sf A}|\hbi{\sf B}|\hbi{\sf C}|\hbi{\sf D}}
\def\vjutges{\hbi{\sf A}|\hbi{\sf B}|\hbi{\sf C}|\hbi{\sf D}|\hbi{\sf E}}
\def\viijutges{\hbi{\sf A}|\hbi{\sf B}|\hbi{\sf C}|\hbi{\sf D}|\hbi{\sf E}|\hbi{\sf F}|\hbi{\sf G}}
\def\ixjutges{\hbi{\sf A}|\hbi{\sf B}|\hbi{\sf C}|\hbi{\sf D}|\hbi{\sf E}|\hbi{\sf F}|\hbi{\sf G}|\hbi{\sf H}|\hbi{\sf I}}
\def\xijutges{\hbi{\sf A}|\hbi{\sf B}|\hbi{\sf C}|\hbi{\sf D}|\hbi{\sf E}|\hbi{\sf F}|\hbi{\sf G}|\hbi{\sf H}|\hbi{\sf I}|\hbi{\sf J}|\hbi{\sf K}}
\def\standard{\hbi{W}|\hbi{T}|\hbi{V}|\hbi{F}|\hbi{Q}}
\def\standardUK{\hbi{W}|\hbi{T}|\hbi{F}|\hbi{Q}|\hbi{V}}
\def\llatins{\hbi{C}|\hbi{S}|\hbi{R}|\hbi{P}|\hbi{J}}
\def\llatinsIDSF{\hbi{S}|\hbi{C}|\hbi{R}|\hbi{P}|\hbi{J}}
\def\wtq{\hbi{W}|\hbi{T}|\hbi{Q}}
\def\wtqii{\hbii{W}|\hbii{T}|\hbii{Q}}
\def\iiijutgesii{\hbii{A}|\hbii{B}|\hbii{C}}

\def\rkf#1#2{\hbi{\hbox{\vbox to0pt{\hsize1ex\vss\halign{\hss##\hss\cr
 #2\cr\noalign{\vskip5pt}#1\cr}}}}}

\def\onehalf{\leavevmode\xiipt
 \raise.4ex\hbox{\xpt 1}\hskip-.2ex/\hskip-.15ex\lower.4ex\hbox{\xpt 2}\hskip.2ex}

\def\skl#1{{\,\sf#1\,}}   

\def\qbloc #1#2#3#4#5#6\ska{%
\nspan=#1\multiply\nspan by2\advance\nspan by-1%
\vrule width.2pt%
\hbox{\vbox{\offinterlineskip\halign{%
##&\vrule## width.2pt%
&&\htskip#2{#4##\/}#3\htskip &\vrule## width.2pt\cr
\noalign{\hrule height.8pt}
\Strut&&#5&\cr
\noalign{\hrule height.8pt}
\Strut&&\ignorespaces#6&\cr
\noalign{\hrule height.8pt}
}}}%
\vrule width.2pt}

\def\qnumx#1\ska{\qbloc1\hfil\hfil\sf{$x$}#1\ska} 
\def\qmulticol(#1)#2:#3\ska{\qbloc{#1}\hfil\hfil\rm{\multispan\nspan\it#2\/\hfil}{#3}\ska}

\def\xnumx#1\ska{\bloc1\hfil\hfil\sf{\lower1pt\hbox{\Strut}}{$x$}#1\ska} 

\hborderstrue

\let\xx\bf
\newcommand\hsk{\hskip5pt}
\newcommand\hssk{\hskip2.5pt}
\newcommand\st{$\ast$}
\newcommand\hph{\hphantom{2}}
\newcommand\hq{\hskip.25em}

\paragraph{2.8}
\textbf{Example.}\hskip.5em
As an illustrative example we will consider the final round of a dance\-sport competition. Specifically, we have chosen the Professional Latin Rising Star section of the 2007 Blackpool Dance Festival (Blackpool, England, 25th May 2007). The data were taken from
\texttt{\small http://www.scrutelle. info/results/estelle/2007/blackpool\underl2007/}.
\ensep
As usual, the final was contested by six couples, whose competition numbers were \skl{3},\skl{4},\skl{31},\skl{122},\skl{264} and \skl{238}. Eleven adjudicators ranked their simultaneous performances in four equivalent dances.

The all-round official result was \skl{3}\better \skl{122}\better \skl{264}\better \skl{4}\better \skl{31}\better \skl{238}.
This result comes from the so-called ``Skating System'', whose name reflects a prior use in figure-skating. The Skating System has a first part which produces a separate result for each dance. This~is done mainly on the basis of the median rank obtained by each couple (by the way, this criterion underlies the ``practical'' method that Condorcet was proposing in 1792/93~\cite[ch.\,8]{mu}).
However, the fine properties of this criterion are lost in the second part of the Skating System, where the~all-round result is obtained by adding up the final ranks obtained in the~different dances.

From the point of view of paired comparisons, it makes sense to base the all-round result on the Llull matrix which collects the 44~rankings produced by the 11~adjudicators over the 4~dances.
As~one can see below, in the present case this matrix exhibits several Condorcet cycles, like for instance \skl{3}\better \skl{4}\better \skl{264}\better \skl{3} and \skl{3}\better \skl{122}\better \skl{264}\better \skl{3}, which means that the competition was closely contested.
In consonance with it, the CLC rates obtained below are quite close to each other,\linebreak[3] particularly for the couples \skl{3},\skl{4},\skl{122} and \skl{264}.
By the way, the CLC result orders the contestants differently than the Borda rule, whose associated ordering is \skl{122}\better \skl{3}\better \skl{4}\better \skl{264}\better \skl{31}\better \skl{238}.




Instead of relative scores and margins, the following tables show their absolute counterparts,
\ie without dividing by the total number of votes. This has the virtue of staying with small integer numbers.

\bgroup\smaller
\initsizesmall

\bigskip
\leftline{\hskip1\parindent
\taula{}%
\numx 3/4/31/122/238/264\ska%
\fitaula
%
\hskip.25\parindent
\taula{Original scores}%
\hmulticol(6)\hskip3.5pt$V_{xy}$:
\hsk\sf 3\hsk|\hsk\sf 4\hsk|\hssk\sf 31\hssk|\sf 122|\sf 238|\sf 264:
\st|\xx23|\xx28|\xx23|\xx28|20/
21|\st|\xx23|20|\xx30|\xx24/
16|21|\st|15|\xx25|18/
21|\xx24|\xx29|\st|\xx28|\xx23/
16|14|19|16|\st|19/
\xx24|20|\xx26|21|\xx25|\st\ska%
\fitaula
%
\hskip.25\parindent
\taula{Indirect scores\hskip-5em}%
\hmulticol(6)\hskip3.5pt\smash{\hbox{$V^\ast_{xy}$}}:
\hsk\sf 3\hsk|\hsk\sf 4\hsk|\hssk\sf 31\hssk|\sf 122|\sf 238|\sf 264:
\st|23|\xx28|23|\xx28|23/
\xx24|\st|\xx24|23|\xx30|\xx24/
21|21|\st|21|\xx25|21/
\xx24|\xx24|\xx29|\st|\xx28|\xx24/
19|19|19|19|\st|19/
\xx24|23|\xx26|23|\xx25|\st\ska%
\fitaula
\hskip.25\parindent
\taula{\hskip-.8emRanks}%
\hmulticol(1):$r_x$: 4/ 2/ 5/ 1/ 6/ 3\ska
\fitaula
}

\bigskip
\leftline{\hskip1\parindent
\taula{}%
\numx 122/4/264/3/31/238\ska%
\fitaula
%
\hskip.25\parindent
\taula{Indirect margins}%
\hmulticol(6)\hskip2pt\smash{\hbox{$M^\mast_{xy}$}}:
\sf 122|\hsk\sf 4\hsk|\sf 264|\hsk\sf 3\hsk|\hssk\sf 31\hssk|\sf 238:
\st|1|1|1|8|9/
\st|\st|1|1|3|11/
\st|\st|\st|1|5|6/
\st|\st|\st|\st|7|9/
\st|\st|\st|\st|\st|6/
\st|\st|\st|\st|\st|\st\ska%
\fitaula
%
\hskip.25\parindent
\taula{Projected margins}%
\hmulticol(6)\hskip2pt\smash{\hbox{$M^\pi_{xy}$}}:
\sf 122|\hsk\sf 4\hsk|\sf 264|\hsk\sf 3\hsk|\hssk\sf 31\hssk|\sf 238:
\st|1|1|1|3|6/
\st|\st|1|1|3|6/
\st|\st|\st|1|3|6/
\st|\st|\st|\st|3|6/
\st|\st|\st|\st|\st|6/
\st|\st|\st|\st|\st|\st\ska%
\fitaula
\hskip.25\parindent
\taula{Rates}%
\hmulticol(1):$R_x$:
3.3636/
3.3864/
3.4091/
3.4318/
3.5682/
3.8409\ska
\fitaula
}

\egroup

\section{Mathematical setting}\label{sec-setting}

We consider a finite set~$\ist$. Its elements represent the options
which are the matter of a vote. The number of elements of $\ist$
is~$N$.

In order to deal with preferences
we must consider (ordered) \dfc{pairs} of elements of $\ist$.
The pair formed by $a$ and $b$, in
this order, will be denoted simply as~$ab$.
The pairs that consist of two copies of the same element, \ie those
of the form~$aa$, are not relevant for our purposes. So, we will
systematically exclude them from our considerations. This will help
towards a more efficient language.
The~set of all proper pairs, \ie the pairs $ab$ with $a\neq  b$,
will be denoted as~$\tie$. 
Unless we say otherwise, from now on \textit{any~statement 
will be understood to imply the assumption that all the pairs that
appear in it are proper pairs, \ie they belong to~$\tie$}.
\label{conveni}

Besides pairs, we will be concerned also with longer sequences $a_0
a_1 \dots a_n$, which will be referred to as {\df paths}.

\paragraph{3.1}
As one could see in the heuristic outline, 
our developments make use of 
the mathematical concept of (binary) \dfd{relation}.
Stating that two elements $a$ and $b$ are in a certain relation~$\rho$ is equivalent to saying that the pair $ab$ is a~member of a certain set~$\rho$.
Because of what has been said,
unless we say otherwise
we will restrict our attention to strict relations,
\ie relations contained in~$\tie$.
Under the convention made above, we can keep the following standard definitions.%
\footnote{
An~equivalent way to put it is that we accept a~relation $\rho\sbseteq\tie$ as total\,/\,transitive, whenever $\rho\cup\{xx\mid x\in\ist\}$
has the corresponding property in the standard sense.%
}
A relation $\rho\sbseteq\tie$ will be called\,:
\ensep
\dfc{asymmetric} when $ab\in\rho$ implies $ba\not\in\rho$;
\ensep
\dfd{total}, or \dfd{complete}, when $ab\not\in\rho$ implies $ba\in\rho$;
\ensep
\dfc{transitive} when the simultaneous occurrence of $ab\in\rho$
and $bc\in\rho$ implies $ac\in\rho$;
\ensep
a \dfc{partial order} when it is 
transitive and asym\-metric;
\ensep
a \dfc{total order} when it is 
transitive, asymmetric and total;
\ensep
a \dfc{total preorder} when it is 
transitive and total.



\bigskip
Given a relation~$\rho$, we will often consider the relation $\hatrho$ that consists of all pairs $ab$ such that $ba\not\in\rho$;\ensep
$\hatrho$~is called the \dfc{codual} of~$\rho$.
The following lemma collects several properties which are immediate consequences of the definitions:

\begin{lemma}\hskip.5em
\label{st:adjoint}

\iim{a}$\hat{\hatrho} = \rho$.

\iim{b}$\rho\sbset\sigma
\ensep\Longleftrightarrow\ensep \hatsigma\sbset\hatrho$.

\iim{c}$\rho \text{ is asymmetric}
\ensep\Longleftrightarrow\ensep \rho\sbseteq\hatrho
\ensep\Longleftrightarrow\ensep \hatrho \text{ is total}$.

\iim{d}$\rho \text{ is total}
\ensep\Longleftrightarrow\ensep \hatrho\sbseteq\rho
\ensep\Longleftrightarrow\ensep \hatrho \text{ is asymmetric}$.
\end{lemma}

\begin{comment}
\noindent
$(\rho\cup\sigma)\,\hat{} = \hatrho \cap \hatsigma$;\quad
$(\rho\cap\sigma)\,\hat{} = \hatrho \cup \hatsigma$;\quad
$\rho\cup\hatrho$ is total;\quad
$\rho\cap\hatrho$ is asymmetric;\hfil\break
$\rho$ transitive $\Rightarrow$ the partial order $\rho\cap\hatrho$ is compatible with the equivalence relation $\rho\cap\rho'$.\par
\end{comment}

\begin{comment}
\noindent $(\rho\cup\sigma)\,\hat{} = \hatrho \cap
\hatsigma$;\quad $(\rho\cap\sigma)\,\hat{} = \hatrho \cup
\hatsigma$;\quad $\rho\cup\hatrho$ is total;\quad
$\rho\cap\hatrho$ is asymmetric;\hfil\break $\rho$ transitive
$\Rightarrow$ the partial order $\rho\cap\hatrho$ is compatible
with the equivalence relation $\rho\cap\rho'$.\par
\end{comment}

\bigskip
The \dfc{transitive closure} of~$\rho$,
which we will denote as $\tcl\rho$, is defined as follows:
$ab\in\tcl\rho$ \ifoi there exists a path $a_0 a_1 \dots a_n$
from $a_0=a$ to $a_n=b$ such that $a_ia_{i+1}\in\rho$ for every $i$.
\ensep
$\tcl\rho$ is the minimum transitive relation that contains $\rho$.
\ensep
The transitive-closure operator is easily seen to~have the following
properties:\ensep
$\tcl\rho\sbseteq\tcl\sigma$ whenever $\rho\sbseteq\sigma$;\ensep
$\tclp{\rho\cap\sigma}\sbseteq(\tcl\rho)\cap(\tcl\sigma)$;\ensep
$(\tcl\rho)\cup(\tcl\sigma)\sbseteq\tclp{\rho\cup\sigma}$;\ensep
$\tclp{\tcl\rho}=\tcl\rho$.



\bigskip
A subset $\cst\sbseteq\ist$ is said to be {\df autonomous} for a
relation~$\rho$ when,\linebreak[3] for any~$x\not\in\cst$, having
$ax\in\rho$ for some $a\in\cst$ implies $bx\in\rho$ for any
$b\in\cst$, and similarly, having $xa\in\rho$ for some $a\in\cst$
implies $xb\in\rho$ for any $b\in\cst$ (see for instance \cite{bm}). \ensep
On~the other hand, $\cst\sbseteq\ist$ will be said to be an {\df interval}
for a relation~$\rho$ when the simultaneous occurrence of
$ax\in\rho$ and $xb\in\rho$ with $a,b\in\cst$ implies $x\in\cst$.
\ensep The following facts are easy consequences of the definitions:
\ensep If~$\rho$ is asymmetric and $\cst$ is autonomous for~$\rho$
then $\cst$ is an interval for~$\rho$. \ensep If~$\rho$ is total and
$\cst$ is an interval for~$\rho$ then $\cst$ is autonomous
for~$\rho$. \ensep As a corollary, if~$\rho$ is total and
asymmetric, then $\cst$ is autonomous for~$\rho$ \ifoi \,it is an
interval for that relation. \ensep Later on we will make use of the
following fact, which is also an easy consequence of the
definitions:

\begin{lemma}\hskip.5em
\label{st:clusterLemma}$\cst$~is autonomous for~$\rho$
\,$\Longleftrightarrow$\, $\cst$~is autonomous for~$\hatrho$.
\end{lemma}

\smallskip
When $\cst$ is an autonomous set for~$\rho$, it is natural to
consider a new set~$\istbis$ and a new relation $\rhobis$ by
proceeding in the following way: \ensep $\istbis$ is obtained from
$\ist$ by replacing the set $\cst$ by a single element~$\clustit$,
\ie $\istbis = (\ist\setminus\cst)\cup\{\clustit\}$; \ensep
for~every $x\in\ist$, let us denote by $\contr{x}$ the element
of~$\istbis$ defined by $\contr{x} = \clustit$ if~$x\in\cst$ and by
$\contr{x} = x$ if~$x\not\in\cst$; \ensep with this notation,
$\rhobis$ is defined by putting $\contr{x}\contr{y}\in\rhobis$ \ifoi
$xy\in\rho$ whenever $\contr{x}\neq\contr{y}$ (this definition is
not ambiguous since $\cst$ is autonomous for~$\rho$). \ensep We~will
refer to this operation as the {\df con\-traction} of $\rho$ by the
autonomous set $\cst$.

\bigskip
Any relation can be interpreted as expressing a system of qualitative preferences:
\ensep
having $xy\in\rho$ and $yx\not\in\rho$ means that $x$ is preferred to~$y$;\ensep
having both $xy\in\rho$ and $yx\in\rho$ means that $x$ is tied with~$y$;\ensep
having neither $xy\in\rho$ nor $yx\in\rho$ means that no information is given about the preference between $x$ and~$y$.

From this point of view, it is quite natural to rank the different $x\in\ist$ 
by~taking into account\, the number of~$y$ such that $x$ is preferred to $y$\,
as~well as the number of~$y$ such that $x$ is tied with $y$.
More precisely, it makes sense to define the \dfc{rank} of $x$ in a relation~$\rho$ by the formula
\begin{equation}
\label{eq:rank}
\rank{x} \,=\, N \,-\, \big|\{\,y\mid xy\cd\in\rho,\ yx\cd{\not\in}\rho\}\big| \,-\, {\textstyle\frac12}\,\,\big|\{\,y\mid xy\cd\in\rho,\ yx\cd\in\rho\}\big|.
\end{equation}
Ranking by $r_x$ is often considered in connection with tournaments. Such a method is known as the Copeland rule (see for instance \cite[p.\,206--209]{t6}, where the tie-splitting term is not present ---since ties are not occurring--- and an equivalent formulation is used).
\ensep
The next lemma is an easy consequence of the definitions.
Its second part justifies using the term ‘rank’.

\begin{lemma}\hskip.5em
\label{st:rank}
If $\rho$~is a partial order, then having $xy\in\rho$ implies $\rank{x}<\rank{y}$.
\ensep
If $\rho$~is a total order, then having $xy\in\rho$ is equivalent to $\rank{x}<\rank{y}$; 
in fact, $\rank{x}$~coincides then with the ordinal number that gives the  position of $x$ in $\rho$.
\end{lemma}

\paragraph{3.2}
The Llull matrix of a vote, as well as the analogous matrices formed respectively by the indirect scores and by the projected scores, are all of them particular instances 
of the abstract notion of \dfc{valued relation}
(also called `fuzzy relation').
In fact, a valued relation on~$\ist$ means simply a mapping~$v$
whereby every pair $xy\in\tie$ is assigned a score $v_{xy}$
in the interval~$[0,1]$.
The~score~$v_{xy}$ measures ``how much'' is $x$ related to~$y$.
Here and in the following we keep the notation and terminology 
introduced in the preceding sections.

\bgroup 
\advance\abovedisplayskip by-.5mm 
\advance\belowdisplayskip by-.5mm 

Most of the notions that are associated with relations
can be generalized to valued relations
(sometimes in several different ways).
Some of these generalized notions have already appeared in the heuristic outline of~\secpar{\ref{sec-heuristic}}. For~instance, we have already remarked that the indirect scores~$\isc_{xy}$ generalize the notion of transitive closure. Another generalized notion, namely that of autonomous set, will appear in~\secpar{\ref{sec-clons}}.


Here we will only remark that for our purposes the valued analogues of asymmetry and totality are respectively the conditions $v_{xy}+v_{yx}\le1$ and $v_{xy}+v_{yx}\ge1$.

\paragraph{3.3}
In general terms, the problem of preference aggregation
deals with\linebreak[3]
valued relations whose scores satisfy the condition $v_{xy}+v_{yx}\le1$.
The set of all such objects will be denoted by $\Omega$.
So, $\Omega = \{\,v\in[0,1]^{\tie} \mid v_{xy}+v_{yx}\le1\,\}$.
The complete case corresponds to the subset $\Gamma$ determined by the equality $v_{xy}+v_{yx}=1$.
So, $\Gamma = \{\,v\in[0,1]^{\tie} \mid v_{xy}+v_{yx}=1\,\}$.


The sets $\Omega$ and $\Gamma$ are respectively the fields of variation of the collective Llull matrix of a vote in the general case and in the complete one. 
On the other hand, the individual votes can also be viewed as belonging to these sets. In~fact, any qualitative expression of preferences (not necesarily transitive) between the elements of $\ist$ can be represented as a Llull matrix, \ie an element of $\Omega$, by putting 
\begin{equation}
\label{eq:binrelmatrix}
v_{xy} =
\begin{cases}
1, &\text{if $x$ is preferred to~$y$},\\
1/2, &\text{if $x$ is tied with~$y$},\\
0, &\text{%
\vtop{\hsize83mm\parindent0pt
\strut if\hskip.5em either $y$ is preferred to~$x$\hskip.5em
or no information is given about the preference between $x$ and $y$.\strut}}
\end{cases}
\end{equation}
Such a mapping satisfies $v_{xy}+v_{yx}=1$ whenever we are in the complete case, \ie when either a preference or a tie is expressed about each pair of options.
\ensep
According to its definition, the collective Llull matrix 
is simply the center of gravity of the~distribution
of individual votes:
\begin{equation}
\label{eq:cog}
v_{xy} = \sum_k \alpha_k\, \vk_{xy},
\end{equation}
where $\vk$ are the Llull matrices associated with the individual votes, and $\alpha_k$ are the corresponding relative frequencies or weights.

\egroup 

\begin{comment}
Notice also that the $\arank{x}$ defined by (\ref{eq:avranksfromscores}) can always be interpreted as $\arank{x} = \sum_k \alpha_k\, \rank{x}^k$, where $\rank{x}^k$ is defined as in \secpar{3.1}. This is true even if the votes are not total orders.
\end{comment}


Since the individual votes play no other role than contributing to the collective Llull matrix as described by the preceding equation,
generally speaking 
there is no need to restrict them to express qualitative preferences only,
but one can allow them to express valued preferences,
\ie to be\linebreak[3] arbitrary elements of~$\Omega$ (of $\Gamma$ in the complete case).
Such a possibility makes sense in that the individual opinions may already be the result of aggregating a variety of criteria.


The preceding idea of representing a system of valued preferences as an element of~$\Omega$ appears already in~\cite{jl}.



\section{The indirect scores and the associated comparison relation}\label{sec-indirect-scores}

Let us recall that the indirect scores~$\isc_{xy}$ are defined in
the following way:
$$
\isc_{xy} \,=\, \max\,\{v_\alpha\mid \text{$\alpha$ is a path \,$x_0
x_1 \dots x_n$\, from $x_0=x$ to $x_n=y$}\,\},
$$
where the score~$v_\alpha$ of a path $\alpha=x_0 x_1 \dots x_n$ is
defined as
$$
v_\alpha \,=\, \min\,\{v_{x_ix_{i+1}}\mid 0 \le i < n\,\}.
$$

\remark
The matrix of indirect scores $\isc$ can be viewed as a power of $v$
(supplemented with $v_{xx} = 1$) for a matrix product defined in the
following way: $(vw)_{xz} = \max_y \min(v_{xy},w_{yz})$. More
precisely, $\isc$ coincides with such a power for any exponent
greater than or equal to~$N-1$.

\vskip2pt 
\medskip
\begin{lemma}\hskip.5em
\label{st:minInequality} The indirect scores satisfy the following
inequalities:
\begin{equation}
\isc_{xz} \ge \min\,(\isc_{xy}, \isc_{yz})\quad\hbox{for any
$x,y,z$.}
\label{eq:minInequality}
\end{equation}
\end{lemma}

\begin{proof}\hskip.5em
Let $\alpha$ be a path from $x$ to $y$ such that $\isc_{xy} =
v_\alpha$; let $\beta$ be a path from $y$ to $z$ such that
$\isc_{yz} = v_\beta$. Consider now their
concatenation~$\alpha\beta$. Since $\alpha\beta$~goes from $x$ to
$z$, one has $\isc_{xz} \ge v_{\alpha\beta}$. On the other hand, the
definition of the score of a path ensures that $v_{\alpha\beta} =
\min\,(v_\alpha, v_\beta)$. Putting these things together gives the
desired result.
\end{proof}

\medskip
\begin{lemma}\hskip.5em
\label{st:indirectEqualToDirect} Assume that the original scores
satisfy the following inequalities:
\begin{equation}
v_{xz} \ge \min\,(v_{xy}, v_{yz})\quad\hbox{for any $x,y,z$.}
\label{eq:minInequalityIndirect}
\end{equation}
In that case, the indirect scores coincide with the original ones.
\end{lemma}

\begin{proof}\hskip.5em
The inequality $\isc_{xz}\ge v_{xz}$ is an immediate consequence of
the definition of $\isc_{xz}$. The converse inequality can be
obtained in the following way: Let $\gamma=x_0x_1x_2 \dots x_n$ be a
path from~$x$ to~$z$ such that $\isc_{xz} = v_\gamma$. By virtue
of~(\ref{eq:minInequalityIndirect}), we have
$$
\min\,\left(\,v_{x_0x_1}, v_{x_1x_2}, v_{x_2x_3}, \dots,
v_{x_{n-1}x_n}\right) \,\,\le\,\, \min\,\left(\,v_{x_0x_2},
v_{x_2x_3}, \dots, v_{x_{n-1}x_n}\right).
$$
So, $\isc_{xz} \le v_{\gamma'}$ where $\gamma'=x_0x_2 \dots x_n$. By
iteration, one eventually gets $\isc_{xz} \le v_{xz}$.
\end{proof}

\remark
On the basis of the preceding results it makes sense to take condition~(\ref{eq:minInequalityIndirect}) as the definition of transitivity for 
a valued relation $(v_{xy})$.\break
\ensep
The matrix of indirect scores $\isc$ can be characterized as the lowest one that lies above $v$ and satisfies such a notion of transitivity; a proof of this fact ---in a more general setting--- will be found in \cite[Theorem~3.3]{dp}. 

\begin{comment}
\noindent Remark: In order to obtain $\isc_{xz}=v_{xz}$ it suffices
that inequality~(\ref{eq:minInequalityIndirect}) holds for all pairs
of the form $xy$ or alternatively for all those of the form $yz$.\par
\end{comment}

\medskip
\begin{theorem}[Schulze, 1998~%
\hbox{\brwrap{\bibref{sc}\,b}}; 
\,see also \hbox{\brwrap{\dbibref{t6}{p.\,228--229}}}%
]
\label{st:transSchulze}
The indirect comparison relation $\icr = \crl(\isc)$ is a partial order.
\end{theorem}

\begin{proof}\hskip.5em
Since $\crl(\isc)$ is clearly asymmetric, it is only a matter of
showing its transitivity.\ensep We will argue by contradiction. Let
us assume that\linebreak $xy\in\crl(\isc)$ and $yz\in\crl(\isc)$, but
$xz\notin\crl(\isc)$. This means respectively that\ensep
(a)~$\isc_{xy} > \isc_{yx}$\ensep and\, (b)~$\isc_{yz} >
\isc_{zy}$,\ensep but\, (c)~$\isc_{zx} \ge \isc_{xz}$. On the other
hand, Lemma~\ref{st:minInequality} ensures also that\ensep
(d)~$\isc_{xz} \ge \min\,(\isc_{xy}, \isc_{yz})$. We~will 
distinguish two cases depending on which of the last two quantities
is smaller:\ensep (i)~$\isc_{yz} \ge \isc_{xy}$;\ensep
(ii)~$\isc_{xy} \ge \isc_{yz}$.

\halfsmallskip Case~(i)\,: $\isc_{yz} \ge \isc_{xy}$.\ensep We will
see that in this case (c) and (d) entail a contradiction with~(a).
In fact, we have the following chain of inequalities: $\isc_{yx} \ge
\min\,(\isc_{yz}, \isc_{zx}) \ge \min\,(\isc_{yz}, \isc_{xz}) \ge
\min\,(\isc_{yz}, \isc_{xy}) = \isc_{xy}$,\ensep where we are using
successively: Lemma~\ref{st:minInequality}, (c), (d) and (i).

\halfsmallskip Case~(ii)\,: $\isc_{xy} \ge \isc_{yz}$.\ensep An
entirely analogous argument shows that in this case (c) and (d)
entail a contradiction with~(b). In fact, we have $\isc_{zy} \ge
\min\,(\isc_{zx}, \isc_{xy}) \ge \min\,(\isc_{xz}, \isc_{xy}) \ge
\min\,(\isc_{yz}, \isc_{xy}) = \isc_{yz}$,\ensep where we are using
successively: Lemma~\ref{st:minInequality}, (c), (d) and (ii).
\end{proof}

\section{Admissible orders}\label{sec-admissible}

Let us recall that an admissible order is a total order $\xi$ such
that $\icr \sbseteq \xi \sbseteq \hat\icr$. Let us recall also that
this definition is redundant since each of the two inclusions
implies the other one. So an admissible order is simply a total
order that extends the partial order~$\icr$. The results that follow 
deal with the existence and efficient finding of such extensions.

\medskip
\begin{theorem}[Szpilrajn, 1930~\cite{sz}]\hskip.5em
\label{st:existenceXiThm} Given a partial order~$\rho$ on a finite
set~$\ist$, one can always find a total order~$\xi$ such that $\rho
\sbseteq \xi \sbseteq \hatrho$. If $\rho$ contains neither~$xy$
nor~$yx$, one can constrain $\xi$ to include the pair $xy$.
\end{theorem}

\vskip-3pt 
\medskip
\begin{proposition}\hskip.5em
\label{st:Copeland}
Let $\rho$ be a partial order.
\ensep
Let $\rank{x}$ denote the rank of $x$ in $\rho$ as defined by equation \textup{(\ref{eq:rank})} of~\,\secpar{3.1}. Any total ordering of the elements of~$\ist$ by non-decreasing values of~$\rank{x}$ is an extension of $\rho$.
\end{proposition}
\begin{proof}\hskip.5em
Let $\xi$ be a total order of~$\ist$ for which $x\mapsto\rank{x}$
does not decrease. This means that $xy\in\xi$ implies $\rank{x}\le
\rank{y}$. Now, the contrapositive of the first statement in
Lemma~\ref{st:rank} ensures that $\rank{x}\le\rank{y}$ implies
$yx\not\in\rho$, \ie $xy\in\hatrho$. So $\xi\subseteq\hatrho$, from
which it follows that also $\rho\sbseteq\xi$.
\end{proof}

\remark
The preceding proposition replaces the problem of finding a total order 
that contains~$\rho$ by the similar problem of finding a total order contained in the total preorder~$\hat\rlrating = \{xy\in\tie\mid \rank{x}\le\rank{y}\}$. However, from a practical point of view the latter is a much easier thing to do, since one is guided by the function $x\mapsto\rank{x}$.

\section{The projection}\label{sec-projection}

Let us recall that our rating method is based upon certain projected scores~$\psc_{xy}$. These quantities are obtained through the corresponding margins~$\pmg_{xy}$ by means of the procedure~(\ref{eq:cprojection1}--\ref{eq:cprojection4}). That procedure makes use of an admissible order $\xi$,
whose existence has been dealt with in the preceding section,
and it assumes $xy\in\xi$.

\medskip
\begin{lemma}\hskip.5em
\label{st:intervals}
The projected margins $\pmg_{xy}$ have the following properties:
\begin{gather}
\hskip-0.8em 
0 \,\le\, \,\pmg_{xy} \,\le\, 1,\qquad
\hskip2.4em 
\hbox{whenever $x\rxi y$}.
\label{eq:intervalsa}
\\[2.5pt]
\begin{repeated}{eq:equaltomax}
\pmg_{xz} \,=\, \max\,(\pmg_{xy}, \pmg_{yz}),\qquad \hbox{whenever $x\rxi y\rxi z$}.
\end{repeated}
\label{eq:intervalsb}
\end{gather}
\end{lemma}
\begin{proof}\hskip.5em
Both properties are immediate consequences of (\ref{eq:cprojection1}--\ref{eq:cprojection4}) and the fact that $0\le\img_{xy}\le1$.
\end{proof}

\medskip
\begin{theorem}\hskip.5em
\label{st:independenceOfXi}
The projected scores do not depend on the admissible
order~$\xi$ used for their calculation, \ie the value of\,
$\psc_{xy}$ is independent of\, $\xi$ \,for every $xy\in\tie$.
On the other hand, the matrix of the projected scores
in an admissible order~$\xi$
is also independent of\, $\xi$;
\ie if $x_i$ denotes the element of rank~$i$ in~$\xi$,
the value of $\psc_{x_ix_j}$ is independent of\, $\xi$
\,for every pair of indices~$i,j$.
\end{theorem}

\vskip-2pt 
\remark
The two statements say different things since the identity of $x_i$ and $x_j$ may depend on the admissible order~$\xi$.

\begin{proof}\hskip.2em 
Let us consider the effect of replacing $\xi$ by another admissible
order~$\xibis$. In the following, the tilde is systematically used
to distinguish between hom\-olo\-gous objects which are associated
respectively with $\xi$~and~$\xibis$; in particular, such a notation
will be used in connection with the labels of the equations which
are formulated in terms of
the assumed admissible order.

With this terminology, we will prove the two following equalities, which amount to the two statements of the theorem. First,
\begin{equation}
\hskip.75em\pmg_{xy} \,=\, \pmgbis_{xy},\hskip.45em\qquad
\hbox to60mm{for any pair\, $xy\,\ (x\neq y)$.\hfil}
\label{eq:ginvariance}
\end{equation}
Secondly, we will see also that
\begin{equation}
\pmg_{x_ix_j} \,=\, \pmgbis_{\tilde x_i\tilde x_j},\qquad
\hbox to60mm{for any pair of indices\, $ij\,\ (i\neq j),$\hfil}
\label{eq:gijinvariance}
\end{equation}
where $x_i$ denotes the element of rank $i$ in $\xi$, and analogously for~$\tilde x_i$ in $\xibis$.
\ensep

\medskip
Now, it is well known that the set of total order extensions of a given partial order is  always connected through transpositions of consecutive elem\-ents (see for instance \cite[p.\,30]{fh}). Therefore,
it suffices to deal
with the case of two admissible orders $\xi$~and~$\xibis$ which
differ from each other
only by the transposition of two consecutive elements. 
So, we will assume that there are two elements $a$ and $b$ such that the only difference between $\xi$~and $\xibis$ is that $\xi$ contains $ab$ whereas $\xibis$~contains~$ba$.
According to the definition of an admissible order, this implies that
$\img_{ab} = \img_{ba} = 0$.

\newcommand\Ant{P}
\newcommand\ant{p}
\newcommand\Pos{Q}
\newcommand\pos{q}

In order to control the effect of the differences between $\xi$ and $\xibis$, we will make use of the following notation:\ensep
$\ant$ will denote the immediate predecessor of~$a$ in $\xi$; in~this connection, any statement about $p$ will be understood to imply the assumption that the set of predecessors  of~$a$ is not empty. Similarly, $\pos$ will denote   the immediate successor of~$b$ in $\xi$; here too, any statement about $\pos$ will be understood to imply the assumption that the set of successors of $b$ is not empty. So, $\xi$ and $\xibis$ contain respectively the paths $\ant ab\pos$ and $\ant ba\pos$. 
Finally, $x'$~means here the immediate successor of $x$ in $\xi$,
which is the same as in $\xibis$ if $x\ne \ant,a,b$.

\medskip
Let us look first at the superdiagonal intermediate projected
margins~$\ppmg_{hh'}$. According to their definition, namely equation~(\ref{eq:cprojection2}),
$\ppmg_{hh'}$ is the minimum of  a certain set of values of
$\img_{xy}$. In a table where $x$ and $y$ are ordered according to
$\xi$, this set is an upper-right rectangle with lower-left vertex
at $hh'$. Using \smash{$\xibis$} instead of $\xi$ amounts to
interchanging two consecutive columns and the corresponding rows of
that table, namely those labeled by $a$ and $b$. In spite of such a
rearrangement, in all cases but one the underlying set from which
the minimum is taken is exactly the same, so the mininum is the
same. The only case where the underlying set is not the same occurs
for~$h=a$ in the order $\xi$, or $h=b$ in the order $\xibis$; but
then the minimum is still the same because the underlying set
includes $\img_{ab} = \img_{ba} = 0$. So,
\begin{align}
\ppmg_{x_ix_{i+1}} \,&=\, \ppmgbis_{\tilde x_i\tilde x_{i+1}},\qquad \hbox{for any $i=1,2,\dots N\!-\!1$,}
\label{eq:mequivariance}\\[2.5pt]\ppmg_{ab} \,&=\, \ppmgbis_{ba} \,=\, 0.
\label{eq:mabba}
\end{align}
On account of the definition of $\pmg_{x_ix_{j}}$ and
$\pmgbis_{\tilde x_i\tilde x_{j}}$,
(\ref{eq:mequivariance}) results in (\ref{eq:gijinvariance}).

\medskip
Finally, let us see that (\ref{eq:ginvariance}) holds too. To this
effect, we begin by noticing that (\ref{eq:mabba}) is   saying  that
\begin{equation}
\pmg_{ab} \,=\, \pmgbis_{ba} \,=\, 0
\label{eq:4gabba}
\end{equation}
Let us consider now the equation $\pmg_{\ant a} = \pmgbis_{\ant b}$,
which is contained in (\ref{eq:gijinvariance}).
On~account of (\ref{eq:intervalsb}), these equalities entail
\begin{equation}
\pmg_{\ant b} \,=\, \pmg_{\ant a} \,=\, \pmgbis_{\ant b} \,=\,
\pmgbis_{\ant a}.
\label{eq:4gcacb}
\end{equation}
By means of an analogous argument, one obtains also that
\begin{equation}
\pmg_{a\pos} \,=\, \pmg_{b\pos} \,=\, \pmgbis_{a\pos} \,=\,
\pmgbis_{b\pos}.
\label{eq:4gbdad}
\end{equation}
On the other hand, (\ref{eq:gijinvariance}) ensures also that
\begin{equation}
\pmg_{xx'} \,=\, \pmgbis_{xx'},\qquad \hbox{whenever \,$x\neq
\ant,a,b$}.
\label{eq:4gxxp}
\end{equation}
Finally, (\ref{eq:intervalsb}) allows to go from (\ref{eq:4gabba}--\ref{eq:4gxxp}) to the desired general equality~(\ref{eq:ginvariance}).
\end{proof}

\vskip-15pt 
\medskip
\begin{theorem}\hskip.5em
\label{st:propertiesOfProjection}
The projected scores and their asssociated margins satisfy the following properties with respect to any admissible order~$\xi$:

\halfsmallskip\noindent
\textup{(a)}~The following inequalities hold whenever $x\rxi y$ and $z\not\in\{x,y\}$:
\begin{alignat}{2}
\label{eq:pvxyinequality}
\psc_{xy} \,&\ge\, \psc_{yx},
&\qquad
\pmg_{xy} \,&\ge\, 0,
\\[2.5pt]
\label{eq:pvinequalities}
\psc_{xz} \,&\ge\, \psc_{yz},
&\qquad
\psc_{zx} \,&\le\, \psc_{zy},
\\[2.5pt]
\label{eq:pminequalities}
\pmg_{xz} \,&\ge\, \pmg_{yz},
&\qquad
\pmg_{zx} \,&\le\, \pmg_{zy},
\end{alignat}

\halfsmallskip\noindent \textup{(b)} If~$\psc_{xy}=\psc_{yx}$, or
equivalently~$\pmg_{xy}=0$,  then
\textup{(\ref{eq:pvinequalities})} and \textup{(\ref{eq:pminequalities})} are
satisfied all of them with an equality sign.
\end{theorem}

\begin{proof}\hskip.5em
Part~(a).\ensep
Let us begin by noticing that (\ref{eq:pvxyinequality}) reduces to (\ref{eq:intervalsa}). Notice also that 
(\ref{eq:pvinequalities}) follows from (\ref{eq:pminequalities}),
and that (\ref{eq:pminequalities}.1) and (\ref{eq:pminequalities}.2)
are equivalent to each other. 
So, it suffices to prove either (\ref{eq:pminequalities}.1) or (\ref{eq:pminequalities}.2).
We will distinguish
three cases, namely:\ensep (i)~$x\rxi y\rxi z$;\ensep (ii)~$z\rxi
x\rxi y$;\ensep (iii)~$x\rxi z\rxi y$.\ensep
In~case~(i), (\ref{eq:pminequalities}.1) follows from (\ref{eq:intervalsb}).
\ensep
In~case~(ii), (\ref{eq:pminequalities}.2) follows from (\ref{eq:intervalsb}).
\ensep
Finally, in case~(iii) it suffices to use (\ref{eq:intervalsa}) to see  that $\pmg_{xz}\ge 0\ge \pmg_{yz}$.




\smallskip
Part~(b).\ensep Similarly to part~(a), it suffices to prove the statement corresponding to (\ref{eq:pminequalities}.1), \ie that $\pmg_{xy}=0$ implies
$\pmg_{xz}=\pmg_{yz}$. This follows immediately from (\ref{eq:intervalsb}) in cases~(i) and (ii). In case~(iii), (\ref{eq:intervalsb}) allows to~derive that $\pmg_{xz}=\pmg_{zy}=0$, and therefore also the equality $\pmg_{xz}=\pmg_{yz}$.
\end{proof}


\medskip
\begin{proposition}\hskip.5em
\label{st:nochange}
Assume
that there exists a total order $\xi$ such that
the original scores and the associated margins satisfy the following conditions:
\begin{alignat}{2}
\label{eq:vxyinequality}
&v_{xy} \ge v_{yx},\text{ i.e.\ } m_{xy} \ge 0,\qquad &&\hbox{whenever $x\rxi y$},
\\[2.5pt]
\label{eq:mequaltomax}
&m_{xz} \,=\, \max\,(m_{xy},m_{yz}),\qquad &&\hbox{whenever $x\rxi y\rxi z$}.
\end{alignat}
In that case, the projected scores coincide with the original ones.
\end{proposition}


\begin{proof}\hskip.5em
We will begin by showing that
\begin{equation}
m_{xz} \ge \min(m_{xy},m_{yz}),\qquad
\hbox{for any $x,y,z$.}
\label{eq:mgreaterthanmin}
\end{equation}
In order to prove  this inequality  we will distinguish six cases depending on the relative position of $x,y,z$ according to $\xi$:
\ensep
(a)~If $x\rxi y\rxi z$, then (\ref{eq:mgreaterthanmin}) is an immediate consequence of (\ref{eq:mequaltomax}).
\ensep
(b)~If $z\rxi y\rxi x$, then (\ref{eq:mequaltomax}) (with $x$ and $z$ interchanged with each other) gives $m_{zx} = \max\,(m_{zy},m_{yx})$, which owing to the antisymmetric character of the margins is equivalent to (\ref{eq:mgreaterthanmin}) with an equality sign.
\ensep
(c)~If $x\rxi z\rxi y$, then
condition~(\ref{eq:vxyinequality})
guarantees that $m_{xz}\ge0\ge m_{yz}=\min(m_{xy},m_{yz})$.
\ensep
(d)~If $z\rxi x\rxi y$, then we have $m_{xz}\ge m_{yz}=\min(m_{xy},m_{yz})$, where the inequality holds because (\ref{eq:mequaltomax}) ensures that $m_{zy}\ge m_{zx}$, and the equality derives from the hypothesis upon $\xi$.
\ensep
(e,f)~The two remaining cases, namely $y\rxi x\rxi z$ and $y\rxi z\rxi x$, are analogous respectively to (c) and (d).

\halfsmallskip
Now, since we are in the complete case, the scores $v_{xy}$ and the margins $m_{xy}$ are related to each other by the monotone increasing transformation $v_{xy}=(1+m_{xy})/2.$ Therefore, the inequality~(\ref{eq:mgreaterthanmin}) on the margins is equivalent to the following one on the scores:
\begin{equation}
v_{xz} \ge \min(v_{xy},v_{yz}),\qquad
\hbox{for any $x,y,z$.}
\label{eq:vgreaterthanmin}
\end{equation}
According to Lemma~\ref{st:indirectEqualToDirect},
this inequality implies that $\isc_{xy}=v_{xy}$ and therefore $\img_{xy}=m_{xy}$.
\ensep
In particular, $\xi$ is ensured to be an admissible order.

\halfsmallskip
Let us now consider any pair $xy$ contained in $\xi$. By   making use of  (\ref{eq:mequaltomax}) we see that $\ppmg_{xy} = \img_{xy} = m_{xy}$.
As a consequence, the equality $\pmg_{xy} = \max\,\{\, \ppmg_{pp'} \;\vert\; x\rxieq p\rxi y\,\}$ becomes $\pmg_{xy} = \max\,\{\, m_{pp'} \;\vert\; x\rxieq p\rxi y\,\}$. From here, a second application of (\ref{eq:mequaltomax}) allows to derive  that $\pmg_{xy} = m_{xy}$, and therefore $\psc_{xy}=v_{xy}$.
%
\end{proof}

\smallskip
Since conditions (\ref{eq:vxyinequality}--\ref{eq:mequaltomax}) of Proposition~\ref{st:nochange} are included among the properties of the projected Llull matrix according to Lemma~\ref{st:intervals}, one can conclude that they fully characterize the projected Llull matrices, and that the operator $(v_{xy})\mapsto(\psc_{xy})$ really deserves being called a projection:
 
\smallskip
\begin{theorem}\hskip.5em
\label{st:projection}
The operator $P:\Gamma\ni(v_{xy})\mapsto(\psc_{xy})\in\Gamma$ is idempotent, \ie $P^2=P$. Its image $P\Gamma$ consists of the complete Llull matrices $(v_{xy})$ that satisfy $(\ref{eq:vxyinequality}\text{--}\ref{eq:mequaltomax})$ for some total order~$\xi$.
\end{theorem}

\section{The rank-like rates}

Let us recall that the rank-like rates $\rlr{x}$ are given by the formula~(\ref{eq:rrates}), or equivalently by (\ref{eq:rratesfrommargins}).
From these formulas one easily checks that they satisfy condition~\llrr.


\medskip
\begin{lemma}\hskip.5em
\label{st:obsRosa}

\iim{a} If $x\rxi y$ in an admissible order $\xi$,
then $\rlr{x}\le\rlr{y}$.

\iim{b} $\rlr{x}=\rlr{y}$ \,\ifoi\, $\psc_{xy}=\psc_{yx}$.

\iim{c} $\rlr{x}\le\rlr{y}$ \,implies\ \,
the inequalities
$(\ref{eq:pvxyinequality}\text{\,--\,}\ref{eq:pminequalities})$.

\iim{d} $\rlr{x}<\rlr{y}$ \,\ifoi\, $\psc_{xy}>\psc_{yx}$.

\iim{e} $\psc_{xy}>\psc_{yx}$ \,implies\ \, $x\rxi y$ in any admissible order $\xi$.
\end{lemma}
\begin{proof}\hskip.5em
Part~(a).\ensep
It is an immediate consequence of formula~(\ref{eq:rrates}) together with the inequalities (\ref{eq:pvxyinequality}) and (\ref{eq:pvinequalities}.1) ensured by Theorem~\ref{st:propertiesOfProjection}.

\halfsmallskip
Part~(b).\ensep
From (\ref{eq:rrates}) it follows that
\begin{equation}
\label{eq:riguals}
\rlr{y} - \rlr{x} \,=\, (\psc_{xy} - \psc_{yx}) \,+\,
\sum_{\latop{\scriptstyle z\neq x}{\scriptstyle z\neq y}}%
\,(\psc_{xz}-\psc_{yz}).
\end{equation}
Let $\xi$ be an admissible order. By symmetry we can assume $xy\in\xi$. As~a consequence, Theorem~\ref{st:propertiesOfProjection} ensures that the terms of (\ref{eq:riguals}) which appear in parentheses are all of them greater than or equal to zero. So the only possibility for their sum to vanish is that each of them vanishes separately, \ie $\psc_{xy}=\psc_{yx}$ and $\psc_{xz}=\psc_{yz}$ for any $z\not\in\{x,y\}$. Finally, part~(b) of Theorem~\ref{st:propertiesOfProjection}  ensures that all of these equalities hold as soon as the first one is satisfied.

\halfsmallskip
Part~(c).\ensep
When the hypothesis is satisfied as a strict inequality, the result follows by combining the contrapositive of~(a) with part~(a) of Theorem~\ref{st:propertiesOfProjection}.
\ensep
In the case of equality, it suffices to combine (b) with part~(b) of that theorem.

\halfsmallskip
Part~(d).\ensep It follows from (c) and its contrapositive on account of (b).

\halfsmallskip
Part~(e).\ensep It follows from (d) and the contrapositive of (a).
\end{proof}

\medskip
The next theorem characterizes the preference relation determined by the rank-like rates
in terms of the indirect comparison relation~$\icr$ defined in~\secpar{2.4}:
\begin{theorem}[\footnote{We thank an anonymous reviewer for certain remarks that led to the present version of this theorem, which is stronger than the original one.}]\hskip.5em
\label{st:RvsNu}
The rank-like rating given by~\textup{(\ref{eq:rrates})}
is related to the indirect comparison relation~$\icr = \crl(\isc)$
in the following way:
\begin{gather}
\label{eq:referee1}
\rlr{x} < \rlr{y} \ensep\Longleftrightarrow\ensep yx\not\in(\hat\icr)^*,
\\[2.5pt]
\label{eq:referee1bis}
\rlr{x} \le \rlr{y} \ensep\Longleftrightarrow\ensep xy\in(\hat\icr)^*.
\end{gather}
\end{theorem}

\begin{proof}\hskip.5em
The statements (\ref{eq:referee1}) and (\ref{eq:referee1bis}) are equivalent to each other (via the contrapositive of each implication plus a swap between $x$ and $y$). So it suffices to prove (\ref{eq:referee1}). On the other hand, to establish the latter it suffices to prove the two following statements:
\begin{gather}
\label{eq:referee1a}
xy\in(\hat\icr)^* \ensep\Longrightarrow\ensep \rlr{x} \le \rlr{y},
\\[2.5pt]
\label{eq:referee1b}
yx\not\in(\hat\icr)^* \ensep\Longrightarrow\ensep \rlr{x} < \rlr{y}.
\end{gather}

\halfsmallskip
Proof of~(\ref{eq:referee1a}).\ensep
By transitivity, it suffices to consider the case $xy\in\hat\icr$.\ensep
Now, from Theorem~\ref{st:existenceXiThm} one easily sees that $xy\in\hat\icr$ implies that $xy$ belongs to some admissible order~$\xi$. The conclusion that $\rlr{x}\le\rlr{y}$ is then ensured by Lemma~\ref{st:obsRosa}.(a).

\halfsmallskip
Proof of~(\ref{eq:referee1b}).\ensep
Since $(\hat\icr)^*$ is complete, $yx\not\in(\hat\icr)^*$ implies $xy\in(\hat\icr)^*$ and therefore, according to~(\ref{eq:referee1a}), $\rlr{x}\le\rlr{y}$.
\ensep
So, (\ref{eq:referee1b}) will follow if we show that
\begin{equation}
\label{eq:referee1c}
\rlr{x}=\rlr{y} \ensep\Longrightarrow\ensep yx\in(\hat\icr)^*\qquad \hbox{whenever $x\neq y$}.
\end{equation}
Let $\xi$ be an admissible order. Since $\xi\subseteq\hat\icr$, the right-hand side of (\ref{eq:referee1c}) is automatically true if $yx\in\xi$;
so, it remains to consider the case where $xy\in\xi$.\ensep
Let us begin by assuming that $y=x'$.
According to Lemma~\ref{st:obsRosa}.(b), the equality $\rlr{x}=\rlr{x'}$
implies that $\pmg_{xx'}=0$, that is $\ppmg_{xx'}=0$, which means that
there exist $a,b$ such that $a\rxieq x \rxi b$ and $\img_{ab}=0$. Now, the latter
implies that $ba\in\hat\icr$, which can be combined with the fact that $x'b,ax\in\xi\subseteq\hat\icr$ to~derive that $x'x\in(\hat\icr)^*$ (with the obvious
modifications if $x\!=\!a$ or $x'\!=\!b$).\linebreak\ensep
Finally, if we only know that $xy\in\xi$, we can use Lemma~\ref{st:obsRosa}.(a) to see that the equality~$\rlr{x}=\rlr{y}$ implies $\rlr{p}=\rlr{p'}$ for any $p$ such that $x\rxieq p \rxi y$, which reduces the problem to the preceding case.
\end{proof}

\medskip
\begin{corollary}\hskip.5em
\label{st:RvsNuCor}

\iim{a}$\rlr{x} < \rlr{y} \,\Rightarrow\, xy\in\icr$.


\iim{b}If $\hat\icr$ is transitive (which is ensured whenever $\icr$ is total),\hfil\break
then\, $\rlr{x}<\rlr{y} \,\Leftrightarrow\, xy\in\icr$.

\iim{c}If $\icr$ contains a set of the form $\xst\times\yst$ with $\xst\cup\yst=\ist$,\hfil\break
then\, $\rlr{x} < \rlr{y}$ \,for any $x\in\xst$ and $y\in\yst$.
\end{corollary}

\begin{proof}\hskip.5em
Part~(a).\ensep This is an immediate consequence of~(\ref{eq:referee1}) since\linebreak
$yx\not\in(\hat\icr)^* \,\Rightarrow\, yx\not\in\hat\icr \,\Leftrightarrow\, xy\in\icr$.


\halfsmallskip Part~(b).\ensep  It is just a matter of noticing that under the hypothesis that $\hat\icr$ is transitive the right-hand side of (\ref{eq:referee1}) reduces to $xy\in\icr$.

\halfsmallskip Part~(c).\ensep  Let $x\in\xst$ and $y\in\yst$. Since
$\xst\times\yst \sbset \icr \sbseteq \hat\icr$, part~(a) ensures that $\rlr{x}\le\rlr{y}$.
So, it suffices to exclude the possibility that $\rlr{x}=\rlr{y}$. By~using~(\ref{eq:referee1}) one easily sees that this equality would imply $yx\in (\hat\icr)^*$.
In~other words, there would be a path from $y\in\yst$ to $x\in\xst$ entirely contained in $\hat\icr$. Such a path would have to include a pair $ab\in\hat\icr$ with $a\in\yst$ and $b\in\xst$, which is not possible since $ab\in\hat\icr$ means $ba\not\in\icr$.
\end{proof}

\medskip
\begin{proposition}\hskip.5em
\label{st:avranks}
Assume that the votes are total orders.
Assume also that the Llull matrix satisfies the hypothesis of Proposition~\ref{st:nochange}.
In that case, the rank-like rates~$\rlr{x}$ coincide exactly with the mean ranks~$\arank{x}$.
\end{proposition}

\begin{proof}\hskip.5em
Recall that the rank-like rates are related to the projected scores in the same way as the mean ranks are related to the original scores when the votes are total orders~(\secpar{2.5}).\ensep
So, the result follows since Proposition~\ref{st:nochange} ensures that the projected scores coincide with the original ones.
\end{proof}

\section{Continuity}

We claim that the rank-like rates $\rlr{x}$ are continuous functions of the binary scores~$v_{xy}$. The main difficulty in proving this statement lies in the admissible order~$\xi$, which plays a central role in the computations. Since $\xi$ varies in a discrete set, its dependence on the data cannot be continuous at~all. Even so, we claim that the final result is still a~continuous function of the data.

\newcommand\Omegamma{\Gamma}

In this connection, one can consider as data the normalized Llull
matrix~$(v_{xy})$, its domain of variation being the set $\Omegamma$
introduced in \secpar{3.3}.\linebreak 
Alternatively, one can consider
as data the relative frequencies of the possible votes,
\ie the coefficients $\alpha_k$ mentioned also in \secpar{3.3}.

\medskip
\begin{theorem}\hskip.5em
\label{st:continuityThm}
The~projected scores~$\psc_{xy}$ and the~rank-like rates~$\rlr{x}$ depend continuously on the Llull matrix~$(v_{xy})$.
\end{theorem}

\begin{proof}\hskip.5em
The dependence of the rank-like rates on the projected scores is given by formula~(\ref{eq:rrates}), which is not only continuous but
even linear (non-homogeneous).
So we are left with the problem of showing that the projection
$P:(v_{xy})\mapsto(\psc_{xy})$ is continuous. As it has been
mentioned above, this is not so clear, since the projected margins
are the result of certain operations which are based upon an
admissible order~$\xi$  which is determined separately. However, we
will see, on~the one hand, that $P$~is continuous as long as
$\xi$~remains unchanged, and on the other hand, that the results
of~\secpar{\ref{sec-admissible}--\ref{sec-projection}} allow to
conclude that $P$~is continuous on the whole of $\Omegamma$ in spite of
the fact that $\xi$~can change. In the following we will use the
following notation: for~every total order $\xi$, we denote by
$\Omegamma_\xi$ the subset of~$\Omegamma$ which consists of the~Llull
matrices for which $\xi$ is an~admissible order, and we denote by
$P_\xi$ the restriction of $P$ to $\Omegamma_\xi$.

We claim that the mapping $P_\xi$ is continuous for every total order $\xi$.\linebreak 
In~order to check the truth of this statement, one has to go over
the different mappings whose composition defines $P_\xi$
(see~\secpar{2.6}), namely:\ensep
$(v_{xy})\mapsto(\img_{xy})$,\ensep
$(\img_{xy})\mapsto(\ppmg_{xy})$,\ensep
and finally $(\ppmg_{xx'})\mapsto(\pmg_{xy})\mapsto(\psc_{xy})$.\ensep
All of these mappings are certainly continuous since they involve
only additions and substractions as well as the \,$\max$\, and \,$\min$\, operations.

Finally, the continuity of $P$ (and the fact that it is well-defined) is
a~consequence of the following facts (see for instance
\cite[\secpar{2-7}]{mk}):\ensep (a)~$\Omegamma=\bigcup_{\xi}\Omegamma_\xi$; this
is true because of the existence of $\xi$
(Theorem~\ref{st:existenceXiThm}).\ensep (b)~$\Omegamma_\xi$ is a
closed subset of $\Omegamma$; this is true because $\Omegamma_\xi$ is
described by a set of non-strict inequalities which concern
quantities that are continuous functions of $(v_{xy})$ (namely the
inequalities $\img_{xy}\ge 0$ whenever $xy\in\xi$).\ensep (c)~$\xi$
varies over a finite set.\ensep (d)~$P_\xi$ coincides with
$P_{\eta}$ at $\Omegamma_\xi\cap\Omegamma_\eta$, as it is proved in
Theorem~\ref{st:independenceOfXi}.
\end{proof}


\medskip
\begin{corollary}\hskip.5em
\label{st:continuityCor} The rank-like rates depend continuously on
the relative frequency of each possible content of an~individual
vote.
\end{corollary}

\begin{proof}\hskip.5em
It suffices to recall that the Llull matrix $(v_{xy})$ is simply the
center of gravity of the distribution specified by these relative
frequencies (formula~(\ref{eq:cog}) of~\secpar{3.3}).
\end{proof}

\section{Decomposition}

Property~\llrd\ 
is concerned with having a partition of~$\ist$ in two sets $\xst$ and $\yst$ such that each member of $\xst$ is unanimously preferred to any member of $\yst$, that is:
\begin{equation}
\label{eq:v1}
v_{xy} \,=\, 1
\quad \hbox{(and therefore $v_{yx}=0$)}
\quad \hbox{whenever\, $xy\in\xst\times\yst$}.
\end{equation}
According to property~\llrd, to be proved in the present section, in the complete case considered in this article such a situation is characterized by the following equalities:
\begin{alignat}{2}
\rlr{x} \,&=\, \rlrbis{x},\qquad &&\hbox{for all $x\in\xst$},
\label{eq:condrfx}\\[3.5pt]
\rlr{y} \,&=\, \rlrbis{y} \,+\, |X|,\qquad &&\hbox{for all $y\in\yst$},
\label{eq:condrfy}\\
\sum_{x\in\xst}\rlr{x} \,&=\, |\xst|\,(|\xst|+1)/2 &&\null,
\label{eq:condrfs}
\end{alignat}
where $\rlrbis{x}$ and $\rlrbis{y}$ denote the rank-like rates which are determined respectively from the submatrices associated with $\xst$ and $\yst$. More specifically, each of preceding equalities is separately equivalent to (\ref{eq:v1}).

In the following we will continue using a tilde to distinguish
between hom\-olo\-gous objects associated respectively with the
whole matrix and with its submatrices associated with $\xst$ and $\yst$.

\medskip
\begin{lemma}\hskip.5em
\label{lema1-vxy}
Given a partition $\ist=\xst\cup\yst$ in two disjoint nonempty sets, one has the following equivalences:
\begin{equation}
\left.
\begin{array}{c} v_{xy}=1 \\[2.5pt]
\forall\,xy\in\xst\times\yst
\end{array}
\right\} \ \Longleftrightarrow\ \left\{
\begin{array}{c} \img_{xy}\cd=1 \\[2.5pt]
\forall\,xy\in\xst\times\yst
\end{array}
\right\}
\ \Longleftrightarrow\
\left\{
\begin{array}{c}
v_{xy}^\pi=1 \\[2.5pt]
\forall\,xy\in\xst\times\yst
\end{array}
\right.
\label{eq:lema1-vxy}
\end{equation}
\end{lemma}

\begin{proof}\hskip.5em
Assume that $v_{xy}=1$ for all $xy\in\xst\times\yst$. Then $v_{yx}=0$, for all such pairs, which implies that $v_\gamma$ vanishes for any path $\gamma$
which goes from $\yst$ to~$\xst$. This fact, together with the inequality $\isc_{xy}\ge v_{xy}$, entails the
following equalities for all $x\in\xst$ and $y\in\yst$: $\isc_{yx}=0$, $\isc_{xy}=1$, and consequently $\img_{xy}=1$.

Assume now that $\img_{xy}=1$ for all $xy\in\xst\times\yst$. Let
$\xi$ be an admissible order. As an immediate consequence of the
definition, it includes the set $\xst\times\yst$. Let $\last$ be the
last element of $\xst$ according to~$\xi$. From the present
hypothesis it is clear that $\ppmg_{\last\last'}=1$, 
which entails that
$\psc_{xy}=1$  for every $xy\in\xst\times\yst$.

Assume now that $\psc_{xy}=1$ for all $xy\in\xst\times\yst$. Let
$\xi$ be an admissible order. Here too, we are ensured that it
includes the set $\xst\times\yst$; this is so by virtue of
Theorem~\ref{st:propertiesOfProjection}.(a). Let $\last$ be the last
element of $\xst$ according to~$\xi$. From the fact that
$\ppmg_{\last\last'}=\pmg_{\last\last'}=1$, one infers that
$\img_{xy}=1$ for all $xy\in \xst\times\yst$.

Finally, let us assume again that $\img_{xy}=1$ for all
$xy\in\xst\times\yst$. Since\linebreak
$\img_{xy}=\isc_{xy}-\isc_{yx}$ and both terms of this difference
belong to $[0,1]$, the only possibility is $\isc_{xy}=1$ and
$\isc_{yx}=0$, which implies that $v_{yx}=0$.\ensep This equality is
equivalent to $v_{xy}=1$. \ensep
\end{proof}

\medskip
\begin{lemma}\hskip.5em
\label{lema3-vxy}
Condition~\textup{(\ref{eq:v1})} implies, for any admissible order, the following equalities:
\begin{alignat}{2}
\label{eq:ppmgxx}
\ppmg_{xx'} &=\, \ppmgbis_{xx'}, \qquad&&\hbox{whenever $x,x'\in\xst$,}\\[2.5pt]
\label{eq:ppmgyy}
\ppmg_{yy'} &=\, \ppmgbis_{yy'}, \qquad&&\hbox{whenever $y,y'\in\yst$,}
\end{alignat}
\end{lemma}
\begin{proof}\hskip.5em
As we saw in the proof of Lemma~\ref{lema1-vxy}, condition~(\ref{eq:v1}) implies the vanishing of $v_\gamma$ for any path $\gamma$ which goes from $\yst$ to~$\xst$. Besides the conclusions obtained in that lemma, this implies also the following equalities:
\begin{alignat}{3}
\label{eq:iscxx}
\isc_{x\bar x} &= \iscbis_{x\bar x},
\hskip1.75em &\img_{x\bar x} &= \imgbis_{x\bar x},
\hskip1.75em &&\hbox{for all $x,\bar x\in\xst$,}\\[2.5pt]
\label{eq:iscyy}
\isc_{y\bar y} &= \iscbis_{y\bar y},
\hskip1.75em &\img_{y\bar y} &= \imgbis_{y\bar y},
\hskip1.75em &&\hbox{for all $y,\bar y\in\yst$.}
\end{alignat}
Let us fix an admissible order~$\xi$. The second equality of
(\ref{eq:lema1-vxy}) not only ensures that $\xi$ includes the set
$\xst\times\yst$, but it can also be combined with (\ref{eq:iscxx})
and (\ref{eq:iscyy}) to obtain respectively (\ref{eq:ppmgxx}) and
(\ref{eq:ppmgyy}).
\end{proof}

\medskip
\begin{theorem} 
\label{st:condrfPro}
Conditions \textup{(\ref{eq:v1})}, \textup{(\ref{eq:condrfx})}, \textup{(\ref{eq:condrfy})} and \textup{(\ref{eq:condrfs})} are equivalent to each other. 
\end{theorem}

\begin{proof}\hskip.5em
Part~(a):\ensep (\ref{eq:v1}) $\Longrightarrow$
(\ref{eq:condrfx}), (\ref{eq:condrfy}) and (\ref{eq:condrfs}).\ensep
As a consequence of~(\ref{eq:ppmgxx}) and (\ref{eq:ppmgyy}) we get the following equalities:
\begin{alignat}{2}
\label{eq:pmgxx}
\psc_{x\bar x} \,&=\, \pscbis_{x\bar x}, \qquad&&\hbox{for all $x,\bar x\in\xst$,}\\[2.5pt]
\label{eq:pmgyy}
\psc_{y\bar y} \,&=\, \pscbis_{y\bar y}, \qquad&&\hbox{for all $y,\bar y\in\yst$.}
\end{alignat}
On the other hand, Lemma~\ref{lema1-vxy} ensures that
\begin{equation}
\hskip1.7em\psc_{xy} \,=\, 1, \qquad\hskip1.1em\hbox{for all $xy\in\xst\times\yst$.}
\end{equation}
When the projected scores are introduced in (\ref{eq:rrates})
these equalities result in (\ref{eq:condrfx}) and (\ref{eq:condrfy}).
Finally, (\ref{eq:condrfs}) is an immediate consequence of (\ref{eq:condrfx}).

\halfsmallskip
Part~(b):\ensep (\ref{eq:condrfx}) $\Rightarrow$ (\ref{eq:v1});\,  (\ref{eq:condrfy}) $\Rightarrow$ (\ref{eq:v1}).\ensep
On account of formula~(\ref{eq:rrates}), conditions (\ref{eq:condrfx}) and (\ref{eq:condrfy}) are easily seen to be respectively equivalent to the following equalities:
\begin{alignat}{2}
\sum_{\latop{\scriptstyle y\in\ist}{\scriptstyle y\neq x}} \psc_{xy}
\,\,&=\,\,
\sum_{\latop{\scriptstyle \bar x\in\xst}{\scriptstyle \bar x\neq x}}
\pscbis_{x\bar x}
\,+\, |Y|,\qquad &&\hbox{for all $x\in\xst$},
\label{eq:condrfxbis}
\\[2.5pt]
\sum_{\latop{\scriptstyle x\in\ist}{\scriptstyle x\neq y}}\psc_{yx}
\,\,&=\,\,
\sum_{\latop{\scriptstyle \bar y\in\yst}{\scriptstyle \bar y\neq y}}
\pscbis_{y\bar y}
\qquad &&\hbox{for all $y\in\yst$}.
\label{eq:condrfybis}
\end{alignat}
Let us add up respectively the equalities (\ref{eq:condrfxbis}) over $x\in\xst$ and the equalities (\ref{eq:condrfybis}) over $y\in\yst$. Since $\psc_{pq}+\psc_{qp} = \pscbis_{pq}+\pscbis_{qp} = 1$, we obtain
\begin{align}
\sum_{\latop{\scriptstyle x\in\xst}{\scriptstyle y\in\yst}} \psc_{xy}
\,\,&=\,\,
|X|\,|Y|,
\label{eq:condrfxsum}
\\[2.5pt]
\sum_{\latop{\scriptstyle y\in\yst}{\scriptstyle x\in\xst}}\psc_{yx}
\,\,&=\,\,
0.
\label{eq:condrfysum}
\end{align}
Since the projected scores belong to $[0,1]$, the preceding equalities imply respectively
\begin{alignat}{2}
\psc_{xy} \,&=\, 1, \qquad&&\hbox{for all $xy\in\xst\times\yst$,}
\label{eq:condrfxall}
\\[2.5pt]
\psc_{yx} \,&=\, 0, \qquad&&\hbox{for all $xy\in\xst\times\yst$,}
\label{eq:condrfyall}
\end{alignat}
which are equivalent to each other since $\psc_{xy}+\psc_{yx}=1$.
Finally, Lemma~\ref{lema1-vxy} allows to arrive at (\ref{eq:v1}).

\halfsmallskip
\leavevmode
Part~(c):\ensep (\ref{eq:condrfs}) $\Rightarrow$ (\ref{eq:v1}).\ensep
From the definition of $\rlr{x}$ and the fact that $\psc_{xy}\le\psc_{xy}+\psc_{yx}=1$, one easily derives the inequality $\sum_{x\in\xst}\rlr{x} \ge |\xst|\,(|\xst|+1)/2$,\linebreak with equality \ifoi (\ref{eq:condrfxall}) holds. So, the result follows again by\linebreak Lemma~\ref{lema1-vxy}.
\end{proof}


\medskip
\begin{corollary}\hskip.5em
\label{guanyador-perdedor-complet}

\iim{a} $\rlr{x}=1$\,\, \ifoi \,$v_{xy}=1$ for all $y\ne x$.

\iim{b} $\rlr{x}=N$ \ifoi \,$v_{xy}=0$ for all $y\ne x$.
\end{corollary}
\begin{proof}\hskip.5em
It suffices to apply Theorem~\ref{st:condrfPro} to the special cases $X=\{x\}$ and $X=A\setminus\{x\}$. The result can also be obtained directly from Lemma~\ref{lema1-vxy}.
\end{proof}

\renewcommand\uplapar{\vskip-9mm\null}

\uplapar
\section{The Condorcet-Smith principle}\label{sec-condorcet}

\medskip
\begin{theorem}\hskip.5em
\label{st:majPrinciple} Both the indirect majority
relation~$\icr=\crl(\isc)$ and the preference relation 
determined by the rank-like rates comply with the Condorcet-Smith principle:\ensep 
If $\ist$ is partitioned in two sets $\xst$ and $\yst$ with the
property that $v_{xy} > 1/2$ for any $x\in\xst$ and $y\in\yst$,
\ensep
then one has also $xy\in\icr$ and $\rlr{x}<\rlr{y}$
for any such $x$ and $y$.
\end{theorem}

\begin{proof}\hskip.5em
Assume that $x\in\xst$ and $y\in\yst$.
Since $\isc_{xy} \ge v_{xy}$, the hypothesis of the theorem entails that $\isc_{xy} > 1/2$.
\ensep
On the other hand, let $\gamma$ be a path from $y$ to $x$ such that $\isc_{yx} = v_\gamma$; since it goes from $\yst$ to $\xst$, this path must contain at least one link $y_iy_{i+1}$ with $y_i\in\yst$ and $y_{i+1}\in\xst$; now, for this link we have $v_{y_iy_{i+1}} \le 1 - v_{y_{i+1}y_i} < 1/2$, which entails that $\isc_{yx} = v_\gamma < 1/2$.
\ensep
Therefore, we get $\isc_{yx} < 1/2 < \isc_{xy}$, \ie $xy\in\icr$.
\ensep
Finally, the fact that this holds for any $x\in\xst$ and $y\in\yst$ implies, by Corollary~\ref{st:RvsNuCor}.(c), that one has also $\rlr{x}<\rlr{y}$ for any such $x$ and $y$.
\end{proof}



\section{Clone consistency}\label{sec-clons}

Clone consistency (also known as independence of clones) refers to the effect of adding or deleting similar options.
For many voting methods, this may change the outcome in a substantial way.
For instance,
replacing a single option~$c$ by a~set $\cst$ of several options similar to~$c$
may change the result from $c$~being the~winner
to giving the victory to some option outside~$\cst$.
This does not seem right: if~$c$ deserves being chosen when going alone,
then in the second situation the right choice should be
some member of $\cst$.

The notion of similarity that is relevant here can be formalized by the 
concept of autonomous set that was introduced in~\secpar{\ref{sec-setting}.1}.
Recall that $\cst$ being autonomous for a given binary relation means that
each element from outside~$\cst$ relates to all elements of~$\cst$ in the same way.
In the context of voting theory, autonomous sets are often called sets of clones.
\ensep
So, it makes sense to ask for the following property, which we call clone consistency: 
If a set of options is autonomous for each of the individual votes, then:
(a)~this set is also autonomous for the social ranking;\, and\,
(b)~contracting it to a single option in all of the individual votes has no other effect in the social ranking than getting the same contraction.

This requirement was introduced in 1986--87 by Thomas M.~Zavist and T.~Nicolaus Tideman, who also devised a method that satisfies it,
namely the rule of ranked pairs~\cite{ti,zt}.

This section is aimed at proving this property for both the indirect comparison relation $\icr$
as well as the preference relation determined by the rank-like rates.
The core results were obtained by Markus Schulze \hbox{\brwrap{\bibref{sc}\,c\refco
\bibref{scbis}}}.

\paragraph{11.1}
In order to prepare the ground, we need to deal first with certain generalities. To begin with, the notion of an autonomous set will be extended to apply not only to a relation, as defined in~\secpar{\ref{sec-setting}.1}, but also to any valued relation~$(v_{xy})$: A subset $\cst\sbseteq\ist$ will be said to be autonomous for $(v_{xy})$ when
\begin{equation}
v_{ax} = v_{bx},\quad v_{xa} = v_{xb},\qquad \text{whenever
$a,b\in\cst$ and $x\not\in\cst$.}
\end{equation}
This definition can be viewed as an extension of that given in~\secpar{\ref{sec-setting}.1} because of the following fact, which follows easily from the definitions:

\smallskip
\begin{lemma}\hskip.5em
\label{st:clons0} Given a binary relation $\rho$, let $u_{xy}$ and
$v_{xy}$ be the binary scores defined respectively by
\begin{equation}
\label{eq:binrelmatrices} \hskip-.325\textwidth
\vbox{\hsize.3\textwidth$ u_{xy} =
\begin{cases}
1, &\text{if } xy\in\rho,\\
0, &\text{if } xy\notin\rho;
\end{cases}
$}\qquad\vbox{\hsize.3\textwidth$
v_{xy} =
\begin{cases}
1, &\text{if } xy\in\rho \text{ and } yx\notin\rho,\\
1/2, &\text{if }  xy\in\rho \text{ and } yx\in\rho,\\
0, &\text{if } xy\notin\rho.
\end{cases}
$}
\end{equation}
One has the following equivalences:

\iim{a}$\cst$ is autonomous for $\rho$ \ifoi \,$\cst$~is autonomous
for $(u_{xy})$.

\iim{b}$\cst$ is autonomous for $\rho$ \ifoi \,$\cst$~is autonomous
for $(v_{xy})$.
\end{lemma}


\medskip
\begin{lemma}\hskip.5em
\label{st:clons2} Assume that $\cst\sbset\ist$ is autonomous for~$(v_{xy})$. Assume also that either $x$ or $y$, or both, lie outside $\cst$. In this case
\begin{equation*}
\isc_{xy} \,=\, \max\,\{\,v_\gamma\mid\gamma\text{ contains no more
than one element of }\cst\,\}
\end{equation*}
\end{lemma}
\begin{proof}\hskip.5em
It suffices to see that any path $\gamma=x_0 \dots x_n$ from $x_0=x$
to~$x_n=y$ which contains more than one element of~$\cst$ can be
replaced by~another one~$\gammabis$ which contains only one such
element and satisfies $v_{\gammabis}\ge v_\gamma$. Consider first
the case where $x,y\not\in\cst$. In this case it will suffice to
take $\gammabis=x_0\dots x_{j-1}x_k\dots x_n$, where $j =
\min\,\{\,i\mid x_i\in\cst\,\}$ and $k = \max\,\{\,i\mid
x_i\in\cst\,\}$, which obviously satisfy $0<j<k<n$. Since
$x_{j-1}\not\in\cst$ and $x_j,x_k\in\cst$, we~have
$v_{x_{j-1}x_j}=v_{x_{j-1}x_k}$, so that
\begin{equation*}
\begin{split}
v_\gamma
\,&=\, \min\,\left(v_{x_0x_1},\dots,v_{x_{n-1}x_n}\right)\\
\,&\le\, \min\,\left(v_{x_0x_1},\dots,v_{x_{j-1}x_j},
v_{x_k x_{k+1}},\dots,v_{x_{n-1}x_n}\right)\\
\,&=\, \min\,\left(v_{x_0x_1},\dots,v_{x_{j-1}x_k}, v_{x_k
x_{k+1}},\dots,v_{x_{n-1}x_n}\right)
\,=\, v_{\gammabis}.
\end{split}
\end{equation*}
The case where $x\not\in\cst$ but $y\in\cst$ can be dealt with in a
similar way by~taking $\gammabis=x_0\dots x_{j-1}x_n$, and
analogously, in the case where $x\in\cst$ and $y\not\in\cst$
it~suffices to take $\gammabis=x_0x_{k+1}\dots x_n$.
\end{proof}

\medskip
\begin{proposition}\hskip.5em
\label{st:clons3}
If~$\cst\sbset\ist$ is autonomous for the scores~$(v_{xy})$,
then $\cst$~is autonomous also for the indirect scores $(\isc_{xy})$.
\end{proposition}
\begin{proof}\hskip.5em
Consider $a,b\in\cst$ and $x\not\in\cst$. Let $\gamma=x_0x_1x_2
\dots x_n$ be a path from $a$ to $x$ such that $\isc_{ax}=v_\gamma$.
By Lemma~\ref{st:clons2}, we can assume that $a$ is the only element
of $\gamma$ that belongs to $\cst$. In particular, $x_1\not\in\cst$,
so that $v_{ax_1} = v_{bx_1}$, which allows to write
\begin{equation*}
\begin{split}
\isc_{ax} \,=\, v_\gamma
\,&=\, \min\,\left(v_{ax_1},v_{x_1x_2},\dots,v_{x_{n-1}x}\right)\\
\,&=\, \min\,\left(v_{bx_1},v_{x_1x_2},\dots,v_{x_{n-1}x}\right)
\,\le\, \isc_{bx}.
\end{split}
\end{equation*}
By~interchanging $a$ and $b$, one gets the reverse inequality
$\isc_{bx}\le\isc_{ax}$ and there\-fore the equality
$\isc_{ax}\cd=\isc_{bx}$. An analogous argument shows that
$\isc_{xa}\cd=\isc_{xb}$.
\end{proof}

\medskip
\begin{corollary}\hskip.5em
\label{st:clonsCor} If~$\cst\sbset\ist$ is autonomous for a relation~$\rho$, then $\cst$~is autonomous also for the transitive closure $\tcl\rho$.
\end{corollary}
\begin{proof}\hskip.5em
Because of Proposition~\ref{st:clons3} and
Lemma~\ref{st:clons0}.(a).
\end{proof}

\paragraph{11.2}
The next results assume that $\cst\sbset\ist$ is autonomous for the
Llull matrix of a vote. Obviously, this assumption is satisfied
whenever $\cst$ is autonomous for all of the individual votes (which
can be allowed to be arbitrary elements of $\Omegamma$ as mentioned 
in~\secpar{3.3}).

\medskip
\begin{theorem}\hskip.5em
\label{st:clons4}
Assume that $\cst\sbset\ist$ is autonomous for the
Llull matrix $(v_{xy})$. Then $\cst$ is autonomous for
the indirect comparison relation $\icr=\crl(\isc)$ as well as for the
total preorder determined by the rank-like rates \textup{(}\ie for the
relation $\hat\rlrating = \{xy\in\tie\mid \rlr{x}\le\rlr{y}\}$\textup{)}.
\end{theorem}
\begin{proof}\hskip.5em
Proposition~\ref{st:clons3} ensures that $\cst$ is autonomous for
the indirect scores $(\isc_{xy})$, from which one easily derives
that $\cst$ is autonomous for the relation~$\icr$. Now, according to
Theorem~\ref{st:RvsNu}, $\hat\rlrating=(\hat\icr)^*$.
So the statement about~$\hat\rlrating$
follows by virtue of Lemma~\ref{st:clusterLemma} and
Corollary~\ref{st:clonsCor}.
\end{proof}

\medskip
\begin{theorem}\hskip.5em
\label{st:clonsp}
Assume that $\cst\sbset\ist$ is autonomous for the
Llull matrix $(v_{xy})$. Then $\cst$ is autonomous also for
the projected scores $(\psc_{xy})$.
\end{theorem}
\begin{proof}\hskip.5em
Since $\psc_{xy}=(1+\pmg_{xy})/2$, it suffices to show that $\cst$ is autonomous for the projected margins $(\pmg_{xy})$.
By Theorem~\ref{st:clons3}, we know that $\cst$ is autonomous for the indirect scores $(\isc_{xy})$, which immediately implies its being autonomous also for the indirect margins $(\img_{xy})$.
So the problem lies at 
showing that the autonomy of $\cst$ is maintained when going
from $(\img_{xy})$ to~$(\pmg_{xy})$ via the procedure (\ref{eq:cprojection1}--\ref{eq:cprojection4}). In the following we let $\xi$ be an admissible order and we distinguish two cases depending on whether $\cst$ is or not an interval for $\xi$.

\halfsmallskip
Assume first that $\cst$ is an interval of $\xi$. Since margins are antisymmetric, in order to prove that $\cst$ is autonomous for $(\pmg_{xy})$ it suffices to show that
\begin{alignat}{2}
\pmg_{xa} &\,=\, \pmg_{xb},\qquad
&&\hbox{for any \,$a,b\in\cst$\, and \,$x\rxi\cst$,}
\label{eq:clonspe}
\\[2.5pt]
\pmg_{ay} &\,=\, \pmg_{by},\qquad
&&\hbox{for any \,$a,b\in\cst$\, and \,$\cst\rxi y$.}
\label{eq:clonspd}
\end{alignat}
In the following we prove~(\ref{eq:clonspe}), the proof of~(\ref{eq:clonspd}) being entirely analogous. If there are no $x\rxi\cst$ there is nothing to prove. Otherwise, let $k$ be the immediate predecessor of the first element of $\cst$. By using (\ref{eq:cprojection3}), one easily sees that (\ref{eq:clonspe}) will follow if we show that
\begin{equation}
\label{eq:clonspsigma}
\ppmg_{cc'} \,\le\, \ppmg_{kk'},\qquad \hbox{for any \,$c\in\cst$.}
\end{equation}
Now, this inequality holds because
\begin{equation*}
\begin{split}
\ppmg_{cc'}
\,=\, \min\,\{\, \img_{pq} \;\vert\; p\rxieq c,\; c'\rxieq q\,\}
\,&\le\, \min\,\{\, \img_{pq} \;\vert\; p\rxieq k,\; c'\rxieq q\,\}\\
\,&=\, \min\,\{\, \img_{pq} \;\vert\; p\rxieq k,\; k'\rxieq q\,\}
\,=\, \ppmg_{kk'},
\end{split}
\end{equation*}
where the inequality is due to the fact that we pass to a smaller set, and the equality that starts the second line holds because
$k'\rxieq q\rxieq c$ implies $q\in\cst$, whereas $p\rxieq k$ implies $p\not\in\cst$, so that $\img_{pq}=\img_{pk'}$ for such $p$ and $q$.

\halfsmallskip
Assume now that $\cst$ is not an interval of $\xi$. That is, there exist $a,b\in\cst$ and $x\not\in\cst$ such that $a\rxi x\rxi b$. This implies that $ax,xb\in\hat\icr$. Since we know that $\cst$ is autonomous for $\icr$, it follows that $cx,xc\in\hat\icr$ for all $c\in\cst$, that is, $\img_{cx}=\img_{xc}=0$ for all $c\in\cst$. From this fact one easily derives, using~(\ref{eq:cprojection2}), that $\ppmg_{pp'}=0$ for all $p\in\ist$ such that $p,p'\in\bar\cst$, where $\bar\cst$ means the minimum interval of $\xi$ that contains $C$. 
Finally, this entails, using (\ref{eq:cprojection3}), that $\pmg_{xy}=0$ for all $x,y\in\bar\cst$, and that any subset of $\bar\cst$, in particular the set $\cst$, is autonomous for $(\pmg_{xy})$.
\end{proof}

\begin{comment}
Pending: General case.
Is $\cst$ autonomous for the relation defined by the fraction-like rates?
\end{comment}

\paragraph{11.3}
Finally, we consider the effect of contracting $\cst$ to a single element.\linebreak 
In this connection we will make use of the notation and definitions
of~\secpar{\ref{sec-setting}.1}, together with the following natural
extension to a system of binary scores: if~$(v_{xy})$ admits $\cst$
as an autonomous set,
the contracted binary scores $(\vbis_{\contr{x}\contr{y}})$ are
characterized by the equality $\vbis_{\contr{x}\contr{y}}=v_{xy}$
whenever $\contr{x}\neq\contr{y}$.
In the following, a tilde is systematically used to distinguish
between hom\-olo\-gous objects associated respectively with
$(\ist,v)$ and~$(\istbis,\vbis)$.

\smallskip
\begin{theorem}\hskip.5em
\label{st:clons5}
Assume that $\cst\sbset\ist$ is autonomous for the
Llull matrix $(v_{xy})$. Then the relation $\icrbis$
coincides with the contraction of $\icr$ by the autonomous
set~$\cst$. Similarly, the relation $\smash{\hat{\rlratingbis}} =
\{xy\in\tiebis\mid \rlrbis{x}\le\rlrbis{y}\}$ coincides with the contraction
of $\hat\rlrating = \{xy\in\tie\mid \rlr{x}\le\rlr{y}\}$ by~$\cst$.
\end{theorem}
\begin{proof}\hskip.5em
We begin by noticing that the operation of taking indirect scores commutes with that of contraction by the autonomous set~$\cst$, \ie $\iscbis_{\contr{x}\contr{y}}=\isc_{xy}$ whenever $\contr{x}\neq\contr{y}$. This is a consequence of Lemma~\ref{st:clons2}. As a consequence, $\icrbis=\crl(\iscbis)$ coincides with the contraction of $\icr=\crl(\isc)$ by $\cst$, and similarly for $\hat{\icrbis}$ and~$\hat\icr$. On~the other hand, a second application of Lemma~\ref{st:clons2}\linebreak 
---to the binary scores associated with the relation $\hat\icr$ by formula~(\ref{eq:binrelmatrices}.1)---\linebreak 
ensures that $(\hat{\icrbis})^*$ is the contraction of $(\hat\icr)^*$.
Finally, in order to see that $\hat{\rlratingbis}$ is the contraction
of~$\hat\rlrating$ it is just a matter of applying Theorem~\ref{st:RvsNu}:
For any $x,y\in A$ such that $\contr{x}\ne\contr{y}$, we have $
\rlrbis{\contr{x}}\le\rlrbis{\contr{y}} \,\Leftrightarrow\,
\contr{x}\contr{y}\in (\hat{\icrbis})^*\,\Leftrightarrow\,
xy\in (\hat\icr)^*\,\Leftrightarrow\, \rlr{x}\le\rlr{y}.$
\end{proof}

\section{About monotonicity}\label{sec-monotonicity}

In this section we consider the effect of raising a particular
option \,$a$\, to a more preferred status in the individual ballots
\textit{without any change in the preferences about the other options}. More
generally, we consider the case where the scores~$v_{xy}$ are
modified into new values $\vbis_{xy}$ such that
\begin{equation}
\label{eq:mona} \vbis_{ay} \ge v_{ay},\quad \vbis_{xa} \le
v_{xa},\quad \vbis_{xy} = v_{xy},\qquad \forall x,y\neq a.
\end{equation}
In such a situation, one would expect the social rates to behave in
the following way, where $y$ is an arbitrary element of
$\ist\setminus\{a\}$:
\begin{gather}
\label{eq:rmona}  \rlrbis{a} \le \rlr{a},
\\[2.5pt]
\label{eq:rrmona}
\rlr{a} < \rlr{y} \,\Longrightarrow\, \rlrbis{a} < \rlrbis{y},\qquad
\rlr{a} \le \rlr{y} \,\Longrightarrow\, \rlrbis{a} \le \rlrbis{y},
\end{gather}
where the tilde indicates the objects associated with the modified
scores. Unfortunately, the rating method proposed in this paper does
not satisfy these conditions, but generally speaking it satisfies
only the following weaker ones:
\begin{gather}
\label{eq:wmona} \rlr{a} < \rlr{y}\ \,\Longrightarrow\,\ \rlrbis{a} \le
\rlrbis{y}.
\\[2.5pt]
\label{eq:smona} (\rlr{a} < \rlr{y},\ \forall y\neq a)\ \,\Longrightarrow\,\
(\rlrbis{a} < \rlrbis{y},\ \forall y\neq a).
\end{gather}
In plain words, (\ref{eq:smona})
is saying that\ensep if $a$ was the only winner for the scores
$v_{xy}$, then it is still the only winner for the scores
$\vbis_{xy}$.



An example exhibiting the lack of property (\ref{eq:rrmona}) can be found in 
\brwrap{
\dbibref{clc}{num\-ber~10 of ``Example inputs''}}.

In this connection, it is interesting to remark that the
method of ranked pairs enjoys also
property (\ref{eq:smona}) \cite[p.\,221--222]{t6},
but it fails at
(\ref{eq:wmona}). A~profile which exhibits
this failure of (\ref{eq:wmona}) 
for the method of ranked pairs is given in 
\brwrap{
\dbibref{clc}{number~9 of ``Example inputs''}}.

\bigskip
The remainder of this section is devoted to 
proving 
properties~(\ref{eq:wmona}) and~(\ref{eq:smona}).

\begin{theorem}\hskip.5em
\label{st:mono} Assume that $(v_{xy})$ and $(\vbis_{xy})$ are
related to each other in accordance with \textup{(\ref{eq:mona})}.
In this case, the following properties are satisfied for any
$x,y\neq a$:
\begin{gather}
\iscbis_{ay} \ge \isc_{ay},\qquad \iscbis_{xa} \le \isc_{xa},
\label{eq:mono1}
\\[2.5pt]
\begin{repeated}{eq:wmona}
\rlr{a} < \rlr{y} \,\Longrightarrow\, \rlrbis{a} \le \rlrbis{y}.
\end{repeated}
\label{eq:mono3}
\end{gather}
\end{theorem}

\begin{proof}\hskip.5em
Let $\icr=\crl(\isc)$ and $\icrbis=\crl(\vbis^\ast)$. Observe
that (\ref{eq:mono1}) implies
\begin{equation}
\imgbis_{ay}\ge \img_{ay},\quad \imgbis_{ya} \le \img_{ya},\quad
\forall y\ne a.\label{eq:mono2}
\end{equation}
We begin by seeing that (\ref{eq:mono3}) will be a consequence of
(\ref{eq:mono1}). In~fact, we have the following chain of implications:
$\rlr{a} < \rlr{y}
\,\Rightarrow\, ay\in \icr
\,\Leftrightarrow\, \img_{ay}>0\linebreak 
\,\Rightarrow\, \imgbis_{ay}>0
\,\Leftrightarrow\, ay\in \icrbis
\,\Rightarrow\, \rlrbis{a} \le \rlrbis{y}$,
where the central implication is provided by~(\ref{eq:mono2}) and the other two strict implications are
guaranteed by Corollary~\ref{st:RvsNuCor}.(a).
So, the problem has been reduced to proving~(\ref{eq:mono1}).

Now, in order to obtain the indirect score for a pair of the form
$ay$ it is useless to consider paths involving $xa$ for some $x\ne
a$, since such paths contain cycles whose omission results in paths
having a larger or equal score. So, the maximum which defines $\isc_{ay}$ is realized by a path which does not involve any pair~$xa$. For such a path $\gamma$ we have
$$\isc_{ay} \,=\, v_\gamma \,\le\, \widetilde v_\gamma \,\le\, \iscbis_{ay},$$
where the first inequality follows directly by (\ref{eq:mona}).
An analogous argument gives $\iscbis_{xa}\le \isc_{xa}$.
\end{proof}

\medskip
\begin{corollary}[\footnote{We thank Markus Schulze for pointing out this fact.}]\unskip
Under the hypothesis of Theorem~\ref{st:mono} one has also the
property \textup{(\ref{eq:smona})}.
\end{corollary}

\begin{proof}\hskip.5em
According to Corollary~\ref{st:RvsNuCor}.(a), the left-hand side of
(\ref{eq:smona}) implies the strict inequality $\isc_{ay} >
\isc_{ya}$ for all $y\neq a$. Now, this inequality can be combined
with (\ref{eq:mono1}) to derive that $\iscbis_{ay} > \iscbis_{ya}$
for all $y\neq a$. Finally, Corollary~\ref{st:RvsNuCor}.(c) with
$\xst=\{a\}$ and $\yst=\ist\setminus\{a\}$ guarantees that the
right-hand side of (\ref{eq:smona}) is satisfied.
\end{proof}

\section{Concluding remarks and open questions}\label{sec-conclusion}

To our knowledge, the existing literature does not offer any other rating method
that combines the quantitative properties of continuity and decomposition
with the Condorcet-Smith principle. As it has been pointed out in the introduction, the latter is quite pertinent when one is interested not only in choosing a winner but in ranking all the alternatives (or in rating them).
\ensep
In particular, the Borda rule (linearly equivalent to rating the options by their mean ranks) satisfies those two quantitative properties, but it does not comply with the Condorcet principle.
\ensep
Quite interestingly, the maximin rates, namely $\sigma_x = \min_{y\ne x}v_{xy}$, combine continuity and a certain form of decomposition with the standard Condorcet principle, that~is, condition~\llmp\ restricted to the case where the subset $\xst$ reduces to a single option.%
\footnote{We thank Salvador Barber\`a for having called our attention to the maximin rating.}
However, they easily fail at this condition when the subset $\xst$ contains several options; a~counterexample is given for instance in \cite[p.\,212--213]{t6}.

\medskip
One may ask whether the CLC~rating method is the only one that satisfies properties 
\llsi--\llrd\ and \llmp.
The answer to this question is surely negative:
Although we are imposing sharp constraints on the rates to be obtained at certain special points of the Llull matrix space, in between these points there is still some degree of freedom (in particular, because rates vary in a continuum).
This remaining freedom might allow for imposing some additional property.
In this connection, we find especially interesting the following

\begin{open}
Can one give a rating method that satisfies \llsi--\llrd\ and \llmp\ together with quantitative monotonicity in the sense that new Llull scores satisfying $(68)$ implies new rates satisfying \hbox{\normalfont(69--70)}?
\end{open}

\medskip
It is also natural to ask whether an alternative rating method with the same properties could be given where the underlying ranking method were not Schulze's method of paths, but some other one. In this connection, we consider especially interesting, because of their properties,
the celebrated rule of Condorcet, Kem\'eny and Young
\hbox{\brwrap{
\bibref{mo}\refsc
\dbibref{t6}{p.\,182--190}%
}},
as well as the ranked-pairs one 
\hbox{\brwrap{
\bibref{ti}\refco
\bibref{zt}\refsc
\dbibref{t6}{p.\,219--223}%
}}:

\begin{open}
Is the rule of Condorcet, Kem\'eny and Young compatible with a con\-tinuous rating method satisfying \llsi--\llrd?
\end{open}

\begin{open}
Is the ranked-pairs rule compatible with a continuous rating method satisfying \llsi--\llrd?
\end{open}

\noindent
On the basis of a previous exploratory work, we believe that the answer to questions~2 and~3 is in both cases negative. More deeply into the question, one can ask:

\begin{open}[\footnote{We thank an anonymous reviewer for having raised this interesting question.}]\hskip.5em
Which properties characterize those ranking methods that can be extended into rating methods satisfying properties~\llsi--\llrd?
\end{open}



\end{document}